\documentclass[oneside,reqno,11pt]{amsart}
\usepackage[utf8]{inputenc}
\usepackage[english]{babel}
\usepackage[T1]{fontenc}
\usepackage{resmes}
\usepackage{microtype}
\usepackage{geometry}
\usepackage{stmaryrd}
\usepackage{amsmath,amssymb,graphicx,amsthm}
\usepackage{mathrsfs}
\usepackage{tikz}
\usetikzlibrary{arrows.meta,calc,patterns}
\usepackage{float}
\usepackage{flafter} 
\usepackage{hyperref}
\hypersetup{
    colorlinks=true,
    linkcolor=blue,
    filecolor=magenta,
    urlcolor=cyan,
}

\theoremstyle{plain}
\newtheorem{theorem}{Theorem}[section]
\newtheorem{lemma}[theorem]{Lemma}
\newtheorem{corollary}[theorem]{Corollary}
\newtheorem{proposition}[theorem]{Proposition}
\theoremstyle{definition}

\newtheorem{hypothesis}[theorem]{Assumption}
\theoremstyle{remark}
\newtheorem{rem}[theorem]{Remark}

\newcommand{\W}{\mathbf{W}}
\newcommand{\cW}{\mathcal{W}}
\newcommand{\cP}{\mathcal{P}}
\newcommand{\R}{\mathbb{R}}

\definecolor{optInk}{HTML}{263746}
\definecolor{optTeal}{HTML}{176F79}
\definecolor{optOrange}{HTML}{B25B23}
\definecolor{optGray}{HTML}{73808A}
\tikzset{
  opt axis/.style={draw=optInk,line width=.65pt,-{Stealth[length=2.1mm]}},
  opt edge/.style={draw=optInk,line width=.75pt},
  opt cut/.style={draw=optTeal,line width=1.1pt},
  opt guide/.style={draw=optGray,line width=.45pt,densely dotted},
  opt leader/.style={draw=optGray,line width=.55pt},
  opt zero cut/.style={draw=optInk,line width=.7pt,dash pattern=on 3pt off 2pt},
  opt dimension/.style={draw=optInk,line width=.5pt,{Stealth[length=1.5mm]}-{Stealth[length=1.5mm]}},
  opt extension/.style={draw=optGray,line width=.35pt},
  opt neighborhood/.style={draw=optOrange,line width=.85pt,dash pattern=on 4pt off 2pt},
  opt point/.style={circle,fill=optInk,inner sep=1.25pt},
}

\newcommand{\OptAxisBreak}[1]{%
  \fill[white] (6.10,-.12) rectangle (6.40,#1);
  \draw[draw=optGray,line width=.45pt,dash pattern=on 2pt off 3pt]
    (6.10,.20)--(6.10,#1);
  \draw[draw=optGray,line width=.45pt,dash pattern=on 2pt off 3pt]
    (6.40,.20)--(6.40,#1);
  \draw[draw=optInk,line width=.65pt] (6.02,-.12)--(6.18,.12);
  \draw[draw=optInk,line width=.65pt] (6.32,-.12)--(6.48,.12);
}

\newcommand{\OptFigureGeometry}{%
  \pgfmathsetmacro{\opteps}{0.34}%
  \pgfmathsetmacro{\opttan}{tan(deg(\opteps))}%
  \pgfmathsetmacro{\optcut}{1-\opteps}%
  \pgfmathsetmacro{\opttop}{\optcut-\opttan/2}%
  \pgfmathsetmacro{\optbottom}{\optcut+\opttan/2}%
}

\begin{document}

\title[Stability of Wasserstein barycenters]{Quantitative stability of classical, entropic, and doubly entropic Wasserstein barycenters}

\author{Guillaume Carlier}
\address{CEREMADE, UMR CNRS 7534, Universit\'e Paris
Dauphine, PSL, Pl. de Lattre de Tassigny, 75775 Paris Cedex 16, France \and MOKAPLAN, INRIA-Paris}
\email{carlier@ceremade.dauphine.fr}

\author{Julien Gu\'erin}
\address{CEREMADE, UMR CNRS 7534, Universit\'e Paris
Dauphine, PSL, Pl. de Lattre de Tassigny, 75775 Paris Cedex 16, France \and MOKAPLAN, INRIA-Paris}
\email{julien.guerin@dauphine.psl.eu}

\author{Quentin M\'erigot}
\address{Université Paris-Saclay, CNRS, Inria,
  Laboratoire de mathématiques d’Orsay, 91405 Orsay, France \and DMA, École normale supérieure, Université PSL, CNRS, 
  75005 Paris, France}
\email{quentin.merigot@universite-paris-saclay.fr}

\date{\today}

\maketitle

\begin{abstract}
In this paper, we show how the recent analysis of (sharp) strong convexity properties of the Kantorovich functional by the third author can be used to establish new global quantitative stability estimates for Wasserstein barycenters and their entropic regularizations with respect to perturbations of the population distribution. For classical Wasserstein barycenters, we improve the state-of-the-art H\"older continuity exponent  with respect to the $\mathcal{W}_1$ metric from $1/6$ to $1/4$, while drastically weakening the assumptions on the population distribution. We also provide an example which shows that the  exponent $1/4$ cannot be improved. Motivated by computational considerations, we also analyze regularizations of Wasserstein barycenters with \emph{outer} and \emph{inner} entropic penalizations. For outer entropic barycenters, we establish novel higher-order stability  results and in particular give a quantitative bound on relative Fisher information. Another contribution of our paper is a strong convexity estimate for an entropic version of the Kantorovich functional via soft Legendre transforms,  enabling us to extend our  stability results to \emph{doubly} entropic barycenters, as introduced by Chizat.

\end{abstract}



\section{Introduction}

\subsection{Wasserstein barycenters}
Since their introduction in \cite{AC11}, barycenters in Wasserstein space have become a widely used tool in applications including image processing \cite{peyretexture}, computer graphics \cite{Solomon}, and statistics \cite{PZ,Sriva,Backhoff2022}. Wasserstein barycenters extend McCann's displacement interpolation \cite{Mccanninterpol} to collections of several, or even infinitely many, probability measures. They therefore provide a natural way of aggregating probability measures while preserving their intrinsic geometry, in contrast with the classical, \emph{flat} notion of an average, which linearly superimposes their densities. 

Throughout the paper, we consider probability measures supported on a fixed closed Euclidean ball $B_R$ of radius $R>0$ in $\R^d$. The space of probability measures $\cP(B_R)$ is endowed with the Wasserstein distance $W_p$, for $p\in\{1,2\}$, both of which metrize weak-$*$ convergence on $\cP(B_R)$. Let $\mathbb{P} \in \cP(\cP(B_R))$ 
be a Borel probability distribution describing the population of measures to be aggregated. A Wasserstein barycenter of $\mathbb{P}$ is defined as a minimizer of the Wasserstein variance functional $F_{\mathbb{P}}$, namely, a solution of
\begin{equation}\label{eq}
\min_{\mu \in \cP(B_R)} F_{\mathbb{P}}(\mu),
\qquad
F_{\mathbb{P}}(\mu)
:=
\int_{\cP(B_R)}
\frac{1}{2} W_2^2(\mu,\rho) \, d\mathbb{P}(\rho).
\end{equation}
Standard compactness arguments ensure the existence of a Wasserstein barycenter, although uniqueness does not always hold. Uniqueness is guaranteed, for instance, when $\mathbb{P}\bigl(\mathcal{P}(B_R) \cap L^1(B_R)\bigr)>0.$
In this case, we denote the unique barycenter by $\mu_{\mathbb{P}}$.

\subsection{Stability of Wasserstein barycenters}
 
Despite their appealing properties, Wasserstein barycenters face an important practical limitation: in most applications, the population distribution $\mathbb{P}$ is either unknown or computationally intractable. It is therefore typically discretized, approximated by a finitely supported measure, or observed only through a finite sample. Questions of consistency, robustness to measurement errors, and convergence of empirical barycenters to their population counterparts therefore raise the following fundamental question: \textit{How quantitatively stable is the solution to Problem \eqref{eq} under perturbations of the population distribution $\mathbb{P}$?}

To formulate this question precisely, the space of population distributions $\cP(\cP(B_R))$ must be equipped with a suitable metric. As is standard in the literature --- see, for instance, \cite{LegouicLoubes,CDM24} --- we use the Wasserstein distance $\cW_1$ associated with the underlying metric space $(\cP(B_R),W_1)$. For $\mathbb{P},\mathbb{Q}\in\cP(\cP(B_R))$, it is defined by
\begin{equation}
\cW_1(\mathbb{P},\mathbb{Q})
:=
\inf_{\gamma\in\Pi(\mathbb{P},\mathbb{Q})}
\int_{\cP(B_R)\times\cP(B_R)}
W_1(\rho,\tilde{\rho})
d\gamma(\rho,\tilde{\rho}),
\end{equation}
where $\Pi(\mathbb{P},\mathbb{Q})$ denotes the set of probability measures on $\cP(B_R)\times\cP(B_R)$ with marginals $\mathbb{P}$ and $\mathbb{Q}$.
Our first objective is to establish quantitative stability estimates under suitable assumptions on $\mathbb{P}$ that, in particular, ensure the uniqueness of its Wasserstein barycenter $\mu_{\mathbb{P}}$. More precisely, we seek constants $C>0$ and $\alpha\in(0,1)$ such that, for every $\mathbb{Q}\in\cP(\cP(B_R))$ and every Wasserstein barycenter $\mu_{\mathbb{Q}}$ of $\mathbb{Q}$,
\[
    W_2(\mu_{\mathbb{P}},\mu_{\mathbb{Q}})
    \leq
    C\,\cW_1(\mathbb{P},\mathbb{Q})^\alpha.
\]

\subsubsection*{Related work}

Wasserstein barycenters are particular instances of Fr\'echet means, which generalize the notion of an average to arbitrary metric spaces through the minimization of an expected squared distance. For qualitative consistency results and laws of large numbers for Wasserstein barycenters, we refer to \cite{LegouicLoubes,BigotKlein}. More general convergence rates for empirical Fr\'echet means were established in \cite{AGP20,LPRS23} by exploiting the Alexandrov geometry of the $2$-Wasserstein space. Applying these results to empirical Wasserstein barycenters, however, requires strong regularity assumptions on $\mu_{\mathbb{P}}$ that are difficult to verify \emph{a priori}. Moreover, these assumptions concern the barycenter itself rather than the population distribution $\mathbb{P}$. By contrast, our setting is deterministic: $\mathbb{Q}$ need not be obtained by sampling from $\mathbb{P}$ and neither distribution is random. 

A distinctive feature of the Wasserstein setting is that the variance functional $F_{\mathbb{P}}$ is convex in the usual, or flat, sense and is therefore amenable to standard tools from convex optimization. Quantitative forms of strong convexity are particularly useful for deriving stability estimates and constitute a central ingredient of our approach. 
To the best of our knowledge, apart from the specific one-dimensional case, the first deterministic quantitative stability estimate for Wasserstein barycenters was obtained in \cite{CDM24}. 
More precisely, it was shown that, if $\mathbb{P}$ assigns positive mass to a family of measures supported on convex sets and whose densities are bounded above and bounded away from zero on their supports, then there exists $C>0$ such that, for every $\mathbb{Q}\in\cP(\cP(B_R))$ and every Wasserstein barycenter $\mu_{\mathbb{Q}}$ of $\mathbb{Q}$,
\begin{equation}\label{eq:holder16}
    W_2(\mu_{\mathbb{P}},\mu_{\mathbb{Q}})
    \leq
    C\,\cW_1(\mathbb{P},\mathbb{Q})^{1/6}.
\end{equation}
We also mention the recent work \cite{han2026quantitative}, which establishes stability of Wasserstein barycenters on Alexandrov spaces with curvature bounded from below.

The proof of \eqref{eq:holder16} in \cite{CDM24} draws heavily on \cite{DM23} and, more broadly, on the rapidly growing literature on the quantitative stability of optimal transport; see \cite{Berman}, \cite{KLM25}, \cite{LetrouitLectureNotes}, and the references therein. A recurring theme in this literature is the quantitative strong convexity of $W_2^2(\cdot,\rho)$, or of an associated dual functional, under suitable assumptions on $\rho$. 

Two recent contributions have substantially advanced this theory. In \cite{cpp26}, Cazelles, Pauwels, and Portales obtained a stability estimate with H\"older exponent $1/5$, which they applied to the quantitative statistical estimation of optimal transport maps. In \cite{Mer26}, the third author established a stability estimate with H\"older exponent $1/4$ for optimal transport maps from a fixed source measure $\rho$ to varying targets, assuming only that $\rho$ has a bounded density. Both works rely on closely related techniques. Moreover, the exponent $1/4$ obtained in \cite{Mer26} is sharp, as shown by a counterexample due to Letrouit \cite{Letrouit26}.  These results have already been applied to the stability of multi-marginal transport from a fixed measure to varying targets \cite{han2026}. This problem is closely related to Wasserstein barycenters, although, unlike the setting considered here, it involves a distinguished reference measure. Building on \cite{Mer26}, we will show that it is possible to improve the exponent in \eqref{eq:holder16} from $1/6$ to $1/4$, while imposing weaker assumptions than those of \cite{CDM24}. We will also provide an example showing that the exponent $1/4$ cannot be improved, at least under our assumptions on the population distribution $\mathbb{P}$.

\subsection{Entropic and doubly entropic Wasserstein barycenters} 

Bigot, Cazelles, and Papadakis observed in \cite{Bigot19} that discretizations of continuous measures may produce Wasserstein barycenters with visible oscillatory artifacts --- which is compatible with their weak-$*$ stability. To rule out such instabilities, they proposed to add an entropic penalty term. Given a regularization parameter $\lambda>0$, the entropic barycenter of $\mathbb{P}$ is defined as the minimizer of
\[
    F_{\mathbb{P},\lambda}
    :=
    F_{\mathbb{P}}+\lambda H,
\]
where $H$ denotes the entropy relative to the Lebesgue measure on $B_R$. The minimizer is always unique, and will be denoted $\mu_{\mathbb{P},\lambda}$. This regularization has several appealing consequences: the barycenter $\mu_{\mathbb{P},\lambda}$ has full support and a Lipschitz log-density. Moreover, as shown in \cite{Bigot19} for finitely supported population distributions, entropic barycenters are stable in relative entropy for every fixed $\lambda>0$. In particular, convergence of the population distributions in $\cW_1$ implies strong $L^1$ convergence of the corresponding densities.  Further characterizations, regularity and stability properties were later established in \cite{CEK21} whose aim was to prove a central limit theorem for these regularized barycenters. 

Although this \emph{outer} entropic regularization yields substantial gains in robustness and regularity, minimizing $F_{\mathbb{P},\lambda}$ remains computationally challenging because it still involves the Wasserstein variance functional $F_{\mathbb{P}}$. To address this difficulty, Chizat introduced in \cite{Chi25} the notion of \emph{doubly entropic} barycenter. In addition to the outer penalty $\lambda H$, the Wasserstein variance is replaced by an average entropic transport cost:
\[
    \mu
    \longmapsto
    \int_{\cP(B_R)}
        \operatorname{OT}_{\lambda}(\mu,\rho)
        \,d\mathbb{P}(\rho),
\]
where $\operatorname{OT}_{\lambda}(\mu,\rho)$ denotes the entropic transport cost between $\mu$ and $\rho$; see Section~\ref{sec:doublyentropicbar} for details. Following \cite{Chi25}, replacing $W_2^2$ by $\operatorname{OT}_{\lambda}$ in the barycenter functional is referred to as \emph{inner} entropic regularization. Although the framework of \cite{Chi25} allows the inner and outer regularizations to have different parameters, we restrict our attention to the case in which the two parameters coincide. 
For fixed positive regularization parameters, \cite{Chi25} also establishes strong stability properties of the doubly entropic barycenter with respect to perturbations of the marginals. Our aim in this work is different, as we wish to establish quantitative stability estimates that remain meaningful as the regularization parameters tend to zero.

\subsection{Contributions}

The goal of this paper is to provide a unified framework, relying on the recent strong convexity estimates from \cite{Mer26}, to establish quantitative stability bounds for  Wasserstein barycenters as well as their entropic and doubly entropic versions.

Our first contribution concerns classical Wasserstein barycenters and improves the H\"older exponent $1/6$ obtained in \cite{CDM24}; see \eqref{eq:holder16}. More precisely, Assumption~\ref{hyp:H0}  only requires that $\mathbb{P}$ gives positive mass to measures admitting an $L^\infty(B_R)$ density, with no positive lower bound or additional regularity requirement. Under this assumption, we prove that
\begin{equation}\label{eq:holder14}
    \forall \mathbb{Q}\in\cP(\cP(B_R)),\qquad W_2(\mu_{\mathbb{P}},\mu_{\mathbb{Q}})
    \leq
    C\,\cW_1(\mathbb{P},\mathbb{Q})^{1/4}.
\end{equation}
We also prove by an explicit counterexample on the barycenter of two very simple measures in dimension $d=2$ that the exponent $1/4$ in \eqref{eq:holder14} is optimal.

Our second contribution concerns entropic barycenters. Under Assumption~\ref{hyp:H0}, we prove
\[
    W_2(\mu_{\mathbb{P},\lambda},\mu_{\mathbb{Q},\lambda})
    \leq
    C\,\cW_1(\mathbb{P},\mathbb{Q})^{1/4},
\]
with the same (independent of $\lambda$) constant $C$ as in inequality \eqref{eq:holder14} for classical barycenters. We then exploit the regularizing effect of the outer entropy to obtain stability estimates in stronger, derivative-sensitive discrepancies. In particular, we prove the following estimate involving the relative Fisher-information 
\begin{equation}\label{eq:stabrelfish}
    \left\|
        \nabla\log(\mu_{\mathbb{P},\lambda})
        -
        \nabla\log(\mu_{\mathbb{Q},\lambda})
    \right\|_{L^2(\mu_{\mathbb{P},\lambda})}
    \leq
    \frac{C}{\lambda}
    \left(
        \cW_1(\mathbb{P},\mathbb{Q})^{1/4}
        +
        \cW_1(\mathbb{P},\mathbb{Q})^{1/16}
    \right),
\end{equation}
as well as the first-order Hellinger-type estimate
\begin{equation}\label{eq:Hellingerorder1}
    \left\|\sqrt{\mu_{\mathbb{P},\lambda}} - \sqrt{\mu_{\mathbb{Q},\lambda}}
    \right\|_{H^1(\mu_{\mathbb{P},\lambda})}
    \leq
    C\left(\frac{\cW_1(\mathbb{P},\mathbb{Q})^{1/4}}{\lambda^{5/4}} +\frac{\cW_1(\mathbb{P},\mathbb{Q})^{1/8}}{\lambda^{9/8}}\right).
\end{equation}
The exponents in these inequalities are probably very far from being sharp. However, they confirm in a quite precise sense the intuition that adding an outer entropy to the barycenter problem enforces a strong  form of robustness.

Our third contribution concerns doubly entropic barycenters. Their analysis requires extending the strong-convexity estimates of \cite{Mer26} to the entropic setting; see Theorems~\ref{th:forte_conv_K_rho_lambda} and \ref{th:forte_conv_OT_lambda_plus_entropie}. These extensions may also be useful in other contexts. As a consequence, if $\mathbb{P}$ satisfies Assumption~\ref{hyp:H0}, we obtain the following combined stability and bias estimate:
\begin{equation}\label{eq:stabdoublyregularizedinformal}
    \forall \mathbb{Q}\in\cP(\cP(B_R)), \qquad W_2(
        \mu_{\mathbb{P},\lambda,\lambda},
        \mu_{\mathbb{Q},\lambda,\lambda}) \leq
    C \left[\cW_1(\mathbb{P},\mathbb{Q}) + \lambda\bigl(|\log\lambda|+1\bigr)
    \right]^{1/4},
\end{equation}
where $\mu_{\mathbb{P}, \lambda, \lambda}$ and $\mu_{\mathbb{Q}, \lambda, \lambda}$ denote doubly entropic barycenters. 
The additional term $\lambda\bigl(|\log\lambda|+1\bigr)$ quantifies the entropic bias and reflects the fact that, unlike $W_2^2$, the entropic transport cost $\operatorname{OT}_\lambda$ generally does not vanish on the diagonal.

The common approach behind these three cases is to exploit coercivity estimates on suitable Bregman divergences associated to the dual problems. In the classical case, we integrate the strong-convexity estimate of \cite{Mer26} over the regular part of the population and combine it with Kantorovich–Rubinstein duality to obtain \eqref{eq:holder14}. Outer entropic regularization preserves this estimate but also comes with an Euler–Lagrange equation that yields higher-order stability \eqref{eq:stabrelfish}--\eqref{eq:Hellingerorder1}. For doubly entropic barycenters, we establish an analogous coercivity estimate for the regularized Kantorovich functional, where the maximum in the definition of the convex conjugate is replaced by a soft-max.

\smallskip

\subsection*{Outline}
The article is organized as follows. Section \ref{sec:classical} is devoted to quantitative stability for classical barycenters, we prove \eqref{eq:holder14} and show by an example that the exponent $1/4$ is sharp in \eqref{eq:holder14}. In section \ref{sec:entropicbar}, we consider entropic barycenters and in particular establish the higher-order stability estimates \eqref{eq:stabrelfish}  and \eqref{eq:Hellingerorder1} for entropic barycenters. Finally, in section \ref{sec:doublyentropicbar}, we study properties of soft Legendre transforms,  establish strong convexity properties of entropic optimal transport costs and deduce from the latter quantitative stability for doubly entropic barycenters i.e. prove \eqref{eq:stabdoublyregularizedinformal}.

\subsection*{Acknowledgments}
QM gratefully acknowledges support from the Agence nationale de la recherche through the PEPR PDE-AI project (ANR-23-PEIA-0004).

\section{Stability of Wasserstein Barycenters}\label{sec:classical}
\subsection{Preliminaries}

Let us start by recalling some classical facts  about the optimal  transport problem (for which we refer to the textbooks \cite{Filippobook} and \cite{VillaniTopics}). 
For $R > 0$, let $B_R$ denote the closed ball of radius $R$ in $\R^d$. Recall that for $p\geq 1$, the $p$-Wasserstein distance, $W_p(\mu, \rho)$, between $\mu\in \cP(B_R)$ and $\rho\in \cP(B_R)$ is defined by 
\[W_p^p(\mu, \rho):=
\inf_{\pi \in \Pi(\mu, \rho)} \int_{B_R\times B_R} |x-y|^p \, d\pi(x,y)\]
where $\Pi(\mu, \rho)$ is the set of transport plans between $\mu$ and $\rho$ i.e. the set of probability measures on $B_R\times B_R$ having $\mu$ and $\rho$ as marginals. In the sequel, we will mainly use the exponents $p=1$ and $p=2$.

Given $\mu$ and $\rho$ two elements of $\cP(B_R)$, Kantorovich duality formula (for $W_2$) reads
\begin{equation}\label{eq:Kantorovichdf}
\frac{1}{2} W_2^2(\mu, \rho)=\sup_{\psi\in C(B_R)} \left\{\int_{B_R} \psi d \mu + \int_{B_R} \psi^c d \rho \right\}
\end{equation}
where $\psi^c$ is the $c$-transform of $\psi$:
\[\psi^c(y):=\min_{x\in B_R} \{ \frac{1}{2}\vert x-y \vert^2-\psi(x)\}.\]
The existence of maximizers $\psi$ (called Kantorovich potentials between $\mu$ and $\rho$) in the formulation \eqref{eq:Kantorovichdf} is well-known and such maximizers will  often be denoted in the sequel (even if not unique) as $\psi_{\mu \to \rho}$. Following \cite{Mer26}, it is convenient to observe that the duality formula \eqref{eq:Kantorovichdf} can also be written in terms of the Kantorovich functional $\mathcal{K}_\rho$:
\[\mathcal{K}_\rho(u) =  \langle u^*, \rho \rangle=\int_{B_R} u^* d \rho\]
where $u^*$ is the Legendre transform of $u\in C(B_R)$ (extended by $+\infty$ outside $B_R$):
\[u^*(y):=\max_{x\in B_R}\left\{ x \cdot y -u(x)\right\}.\]
We refer to \cite{Mer26} for an in-depth study of strong convexity properties of $\mathcal{K}_\rho$ which will be crucial in the sequel. The simple change of variables $\psi\mapsto u=\frac{1}{2} \vert \cdot \vert^2-\psi$ enables to rewrite \eqref{eq:Kantorovichdf} as
\begin{equation}\label{eq:Brenierdf}
\frac{1}{2} W_2^2(\mu, \rho)=\frac{1}{2} m_2(\mu)+ \frac{1}{2} m_2(\rho)-\inf_{u\in C(B_R)} \left\{ \int_{B_R} u d\mu+\mathcal{K}_\rho(u)\right\}
\end{equation}
where $m_2(\mu)$ denotes the second moment of $\mu$:
\[m_2(\mu)=\left\langle |\cdot|^2, \mu \right\rangle.\]
If $\rho$ is absolutely continuous with respect to the Lebesgue measure and $u\in C(B_R)$ then (the convex and Lipschitz function) $u^*$ is differentiable $\rho$ a.e. and $\mathcal{K}_\rho$ is differentiable at $u$ with $\nabla \mathcal{K}_\rho(u)=-\nabla u^*_\# \rho$ (see Proposition 2.2 in \cite{DM23}). Hence for an absolutely continuous $\rho$, setting $u=\frac{1}{2}\vert \cdot \vert^2-\psi_{\mu\to \rho}$, we have $\nabla u^*_\#\rho =\mu$. In other words, the map with a convex potential  $T_{\rho \to \mu}:=\nabla u^*$ is Brenier's  optimal transport (see \cite{Brenier}) between $\rho$ and $\mu$.

Let $\mathbb{P} \in \cP(\cP(B_R))$ be a Borel probability measure on the metric space $(\cP(B_R), W_2)$. 
We associate with $\mathbb{P}$ the functional on $\cP(B_R)$ defined by:
$$ F_{\mathbb{P}}(\mu) = \int_{\cP(B_R)} \frac{1}{2} W_2^2(\mu, \rho) \, d\mathbb{P}(\rho).$$

The variance functional $F_{\mathbb{P}}$ is convex  (for the flat geometry of $\cP(B_R)$) and lower semi-continuous. We define $\mu_{\mathbb{P}}$, a Wasserstein barycenter associated with $\mathbb{P}$, by:
$$ \mu_{\mathbb{P}} \in \mathop{\mathrm{argmin}}_{\mu \in \cP(B_R)} F_{\mathbb{P}}(\mu) = \int_{\cP(B_R)} \frac{1}{2} W_2^2(\mu, \rho) \, d\mathbb{P}(\rho).$$
The existence of a Wasserstein barycenter $\mu_{\mathbb{P}}$ follows immediately from compactness and lower semicontinuity. The uniqueness can be proven under assumptions of the type $\mathbb{P}(\cP(B_R)\cap L^1(B_R)) > 0$ (see \cite{AC11} for weaker assumptions) which guarantee that $F_{\mathbb{P}}$ is strictly convex. Without assumptions on $\mathbb{P}$, it has nevertheless been shown by \cite{CDM24} that one cannot hope for quantitative stability. The uniqueness condition of the barycenters is in fact not even sufficient; requiring $\mathbb{P}$ to assign mass to measures with bounded densities seems to be a good compromise.

\begingroup
\renewcommand{\thetheorem}{($H_0$)}
\begin{hypothesis}\label{hyp:H0}
There exist constants $M_{\mathbb{P}} > 0$ and $\alpha_{\mathbb{P}} > 0$, as well as a Borel set $S \subset \cP(B_R)$, such that $\mathbb{P}(S) \geq \alpha_{\mathbb{P}}$ and for all $\rho \in S$, $\rho$ admits a density satisfying $\|\rho\|_\infty \le M_{\mathbb{P}}$.
\end{hypothesis}
\addtocounter{theorem}{-1}
\endgroup

Assumption \ref{hyp:H0} ensures $\mathbb{P}(L^1(B_R)) > 0$ and thus guarantees uniqueness of the barycenter $\mu_{\mathbb{P}}$.  It is also known (see \cite{AC11}, \cite{CEK21}) that assumption \ref{hyp:H0} ensures that $\mu_{\mathbb{P}}$ is absolutely continuous with a density explicitly bounded by 
\begin{equation}\label{eq:boundclassiccbaryc}
 \|\mu_{\mathbb{P}}\|_\infty \le \frac{M_{\mathbb{P}}}{\alpha^d_{\mathbb{P}}}.
 \end{equation}

\subsection{Stability of Wasserstein barycenters}

Let us now state the main result of this section:

\begin{theorem}[Quantitative stability of the classical barycenter]\label{th:stab_classique}
Let $\mathbb{P}, \mathbb{Q} \in \cP(\cP(B_R))$. Let $\mu_{\mathbb{P}}$ and $\mu_{\mathbb{Q}}$ be Wasserstein barycenters associated with $\mathbb{P}$ and $\mathbb{Q}$, respectively. If $\mathbb{P}$ satisfies \ref{hyp:H0}, then:
\begin{equation}\label{eq:stabbarycclassic}
     W_2(\mu_{\mathbb{P}}, \mu_{\mathbb{Q}}) \le C_{d,R, \alpha_{\mathbb{P}}, M_{\mathbb{P}}} \, \cW_1(\mathbb{P}, \mathbb{Q})^{1/4} 
     \end{equation}
where $C_{d,R, \alpha_{\mathbb{P}}, M_{\mathbb{P}}} = 2^{d+4} \left( \frac{6 \, d^2 \, \beta_d \, M_{\mathbb{P}}}{\alpha_{\mathbb{P}}} \right)^{\frac{1}{4}} R^{\frac{d+3}{4}}$ and $\beta_d$ is the Lebesgue measure of the unit ball.
\end{theorem}

\begin{rem}
This improves the stability result of \cite{CDM24}, which requires  significantly stronger assumptions on $\mathbb{P}$ and where the  exponent is $1/6$ instead of $1/4$.
\end{rem}

The key ingredient of the proof is a quantitative form of the strict convexity of $F_{\mathbb{P}}$, which is obtained by integrating the  estimate very recently obtained  in \cite[Theorem 2.2]{Mer26} and valid for a broad class of probability measures $\rho$. The following strong convexity is a quite direct consequence of \cite[Theorem 2.2]{Mer26}:

\begin{lemma}\label{lem:forte_convexe_W2}
Let $\rho \in \cP(B_R) \cap L^\infty(B_R)$, and $\mu, \mu_0 \in \cP(B_R)$. Let $\psi_0=\psi_{\mu_0\to \rho}$ be a  Kantorovich potential between $\mu_0$ and $\rho$, one has 
\begin{equation}\label{eq:strongconvw21}
    C_{d, R, \|\rho\|_\infty} \left( \frac{1}{2}W_2^2(\rho, \mu) - \frac{1}{2} W_2^2(\rho, \mu_0) - \langle \psi_0, \mu - \mu_0 \rangle\right) \ge  W_2^4(\mu_0, \mu)
\end{equation}
with 
\begin{equation}\label{eq:formulaforC}
 C_{d,R,M}=3 \times 2^{4d+15} R^{d+2} \beta_d d^2 M.
\end{equation}
\end{lemma}

\begin{proof}[Proof of Lemma \ref{lem:forte_convexe_W2}]
Let $\psi$ be a Kantorovich potential between $\mu$ and $\rho$ and let $u, u_0$ denote the Brenier potentials $\frac{1}{2}\vert \cdot \vert^2-\psi$, $\frac{1}{2}\vert \cdot \vert^2-\psi_0$. We thus have  ${\nabla u_0^*}_\# \rho = \mu_0$ and $\nabla u^*_\# \rho = \mu$ and
\begin{align*}
    \frac{1}{2}W_2^2(\rho, \mu) &=\frac{1}{2} m_2(\rho)+ \frac{1}{2} m_2(\mu) - \left\langle  u, \mu \right\rangle - \mathcal{K}_\rho(u) \\
   \frac{1}{2}W_2^2(\rho, \mu_0) & =\frac{1}{2} m_2(\rho)+ \frac{1}{2} m_2(\mu_0) - \left\langle  u_0, \mu_0 \right\rangle - \mathcal{K}_\rho(u_0).
\end{align*}
Hence: 
\begin{align*}
    \frac{1}{2}W_2^2(\rho, \mu) -\frac{1}{2} W_2^2(\rho, \mu_0) - \langle \psi_0, \mu - \mu_0 \rangle 
    &= \mathcal{K}_\rho(u_0) - \mathcal{K}_\rho(u) + \left\langle u_0 - u, \mu \right\rangle \\
    &= \mathcal{K}_\rho(u_0) - \mathcal{K}_\rho(u) + \langle u_0 - u, \nabla u^* \# \rho \rangle \\
    &= \mathcal{K}_\rho(u_0) - \mathcal{K}_\rho(u) - \langle u_0-u, \nabla \mathcal{K}_\rho(u) \rangle \\
    &\ge \frac{1}{ C} \|\nabla u^* - \nabla u_0^*\|_{L^2(\rho)}^4 \ge \frac{1}{ C} W_2^4(\mu, \mu_0),
\end{align*}
where the penultimate inequality stems from the strong convexity estimate for $\mathcal{K}_\rho$ from \cite[Theorem 2.2]{Mer26} (including the explicit value for the constant $C=3 \times 2^{4d+15} R^{d+2} \beta_d d^2 \|\rho\|_\infty$) and the final one from the suboptimality of the transport plan $(\nabla u^*, \nabla u_0^*)_\# \rho$ between $\mu$ and $\mu_0$.
\end{proof}

Now, since $F_{\mathbb{P}} = \int_{\cP(B_R)} \frac{1}{2} W_2^2(\cdot, \rho) \, d\mathbb{P}(\rho)$ is an average of $W_2^2(\cdot, \rho)$, if $\mathbb{P}$ assigns enough mass to good probability measures $\rho \in \cP(B_R) \cap L^\infty(B_R)$, then the strict convexity of $F_{\mathbb{P}}$ can be quantified as follows:

\begin{corollary}\label{cor:forte_conv_FP}
Let $\mathbb{P}$ satisfy assumption \ref{hyp:H0} and let $\mu, \mu_0 \in \cP(B_R)$. Then, for any measurable selection $\rho\mapsto \psi_{\mu_0 \to \rho} \in C(B_R)$ of Kantorovich potentials between $\mu_0$ and $\rho$, defined for $\mathbb{P}$-almost every $\rho$,
$$ F_{\mathbb{P}}(\mu) - F_{\mathbb{P}}(\mu_0) - \int_{\cP(B_R)} \langle \psi_{\mu_0 \to \rho}, \mu - \mu_0 \rangle \, d\mathbb{P}(\rho) \ge \frac{\alpha_{\mathbb{P}}}{C_{d, R, M_{\mathbb{P}}}} \,  \, W_2^4(\mu, \mu_0)$$
where $C_{d, R, M_{\mathbb{P}}}$ is as in \eqref{eq:formulaforC}.
\end{corollary}

\begin{proof}
Integrating the inequality from Lemma \ref{lem:forte_convexe_W2} over all $\rho \in S$ (satisfying $\|\rho\|_\infty \le M_{\mathbb{P}}$):
\[\begin{split}
F_{\mathbb{P}}(\mu)& - F_{\mathbb{P}}(\mu_0) - \int_{\cP(B_R)} \langle \psi_{\mu_0 \to \rho}, \mu - \mu_0 \rangle \, d\mathbb{P}(\rho) \\
&= \int \left[ \frac{1}{2} W_2^2(\rho, \mu) - \frac{1}{2} W_2^2(\rho, \mu_0) - \langle \psi_{\mu_0 \to \rho}, \mu - \mu_0 \rangle \right] d\mathbb{P}(\rho) \\
&\ge \int_S \left[ \frac{1}{2} W_2^2(\rho, \mu) - \frac{1}{2} W_2^2(\rho, \mu_0) - \langle \psi_{\mu_0 \to \rho}, \mu - \mu_0 \rangle \right] d\mathbb{P}(\rho) \\
&\ge \frac{ \alpha_{\mathbb{P}}} {C_{d, R, M_{\mathbb{P}}}} \, W_2^4(\mu, \mu_0).
\end{split}\]
The first inequality uses the non-negativity of the integrand, which relies on the convexity of $W_2^2(\rho,\cdot)$ and the fact that $\psi_{\mu_0\to\rho}$ is a subgradient.
\end{proof}

Let us now proceed to the proof of Theorem \ref{th:stab_classique}.

\begin{proof}[Proof of Theorem \ref{th:stab_classique}]
 Since $\mu_{\mathbb{P}}$ is a minimizer of $F_{\mathbb{P}}$, the first-order condition (see \cite{AC11} for the case of a discrete $\mathbb{P}$ and \cite[Theorem 1.1]{CDM24} for the general case) yields a family $(\psi_{\mu_{\mathbb{P}} \to \rho})_{\rho} \in L^\infty(\mathbb{P}, C(B_R))$ of Kantorovich potentials between $\mu_\mathbb{P}$ and $\rho$ satisfying on $B_R$ the Euler-Lagrange equation:
\begin{equation}\label{eq:eulerlagrangebaryc}
\int_{\cP(B_R)} \psi_{\mu_{\mathbb{P}} \to \rho} \, d\mathbb{P}(\rho) = 0.
\end{equation}
Thus, by Fubini's Theorem and applying Corollary \ref{cor:forte_conv_FP} to $\mu_0 = \mu_{\mathbb{P}}$ and $\mu=\mu_{\mathbb{Q}}$, we get:
$$ F_{\mathbb{P}}(\mu_{\mathbb{Q}}) - F_{\mathbb{P}}(\mu_{\mathbb{P}}) \ge \frac{ \alpha_{\mathbb{P}}} {C_{d, R, M_{\mathbb{P}}}} W_2^4(\mu_{\mathbb{P}}, \mu_{\mathbb{Q}}).$$
By the minimality of $\mu_{\mathbb{Q}}$, we also have:
$$ F_{\mathbb{Q}}(\mu_{\mathbb{P}}) - F_{\mathbb{Q}}(\mu_{\mathbb{Q}}) \ge 0.$$
Thus, summing these inequalities:
\begin{align*}
\frac{ \alpha_{\mathbb{P}}} {C_{d, R, M_{\mathbb{P}}}}  W_2^4(\mu_{\mathbb{P}}, \mu_{\mathbb{Q}}) &\le \left(F_{\mathbb{P}}(\mu_{\mathbb{Q}}) - F_{\mathbb{Q}}(\mu_{\mathbb{Q}})\right) - \left(F_{\mathbb{P}}(\mu_{\mathbb{P}}) - F_{\mathbb{Q}}(\mu_{\mathbb{P}})\right) \\
&= \int_{\cP(B_R)} \left( \frac{1}{2} W_2^2(\mu_{\mathbb{Q}}, \rho) - \frac{1}{2} W_2^2(\mu_{\mathbb{P}}, \rho) \right) d(\mathbb{P} - \mathbb{Q})(\rho) \\
&\le \text{Lip}_{W_1} \left( \frac{1}{2} W_2^2(\mu_{\mathbb{Q}}, \cdot) - \frac{1}{2} W_2^2(\mu_{\mathbb{P}}, \cdot) \right) \cW_1(\mathbb{P}, \mathbb{Q}).
\end{align*}
The last inequality follows from the Kantorovich-Rubinstein inequality. Finally, $\frac{1}{2}W_2^2(\mu, \cdot)$ is $2R$- Lipschitz in $(\cP(B_R), W_1)$. Indeed, if we fix $\psi$ a Kantorovich potential (for the quadratic cost) between $\rho$ and $\mu$, we can assume that $\psi$ is a $c$-transform hence  $2R$-lipschitz. So by Kantorovich duality and Kantorovich-Rubinstein inequality, 
\begin{align*}
 \frac{1}{2}W_2^2(\mu, \rho) - \frac{1}{2}W_2^2(\mu, \tilde{\rho}) &\leq  (\langle \psi^{c}, \mu\rangle  + \langle \psi, \rho\rangle) - (\langle \psi^{c}, \mu\rangle  + \langle \psi, \tilde{\rho}\rangle) \\
&\le  \langle \psi, \rho - \tilde{\rho} \rangle \\
&\le 2R\cdot W_1(\rho, \tilde{\rho}).
\end{align*}
From which we conclude that:
$$ W_2^4(\mu_{\mathbb{P}}, \mu_{\mathbb{Q}}) \le \frac{ 4R C_{d, R, M_{\mathbb{P}}}} {\alpha_{\mathbb{P}}}  \, \cW_1(\mathbb{P}, \mathbb{Q})$$
which, raised to the power $1/4$  exactly yields \eqref{eq:stabbarycclassic} with the constant
\[C_{d,R, \alpha_{\mathbb{P}}, M_{\mathbb{P}}}=\Big(\frac{ 4R C_{d, R, M_{\mathbb{P}}}} {\alpha_{\mathbb{P}}}\Big)^{1/4}
=2^{d+4} \left( \frac{6 \, d^2 \, \beta_d \, M_{\mathbb{P}}}{\alpha_{\mathbb{P}}} \right)^{\frac{1}{4}} R^{\frac{d+3}{4}}. \qedhere \]
\end{proof}

\subsection{Optimality of the exponent in theorem \ref{th:stab_classique}}\label{sec:optimality} 
\begin{figure}[!t]
  \centering
  \resizebox{\linewidth}{!}{
\begin{tikzpicture}[x=1cm,y=1cm,font=\small,text=optInk]
\OptFigureGeometry
\begin{scope}[x=7.7cm,y=7.7cm]
  \fill[optTeal!9] (0,0)--(\optcut,0)--(\optcut,.5)--(\opttop,1)--(0,1)--cycle;
  \fill[optTeal!23] (\optbottom,0)--(1,0)--(1,1)--(\optcut,1)--(\optcut,.5)--cycle;
  \fill[optOrange!14] (\opttop,1)--(\optcut,1)--(\optcut,.5)--cycle;
  \fill[pattern=north east lines,pattern color=optOrange!75]
    (\opttop,1)--(\optcut,1)--(\optcut,.5)--cycle;
  \fill[optOrange!14] (\optcut,0)--(\optbottom,0)--(\optcut,.5)--cycle;
  \fill[pattern=north west lines,pattern color=optOrange!75]
    (\optcut,0)--(\optbottom,0)--(\optcut,.5)--cycle;
  \draw[opt edge] (0,0) rectangle (1,1);
  \draw[opt guide] (0,.5)--(\optcut,.5);
  \draw[opt zero cut] (\optcut,0)--(\optcut,1);
  \draw[opt cut] (\opttop,1)--(\optbottom,0);
  \draw[opt axis] (0,0)--(1.09,0) node[right] {$x_1$};
  \draw[opt axis] (0,0)--(0,1.13) node[above] {$x_2$};
  \node[below left] at (0,0) {$0$};
  \node[left] at (0,1) {$1$};
  \node[left] at (0,.5) {$\tfrac12$};
  \node[below] at (1,0) {$1$};
  \node[below=3pt] at (\optcut,0) {$1-\varepsilon$};
  \node[opt point] at (\opttop,1) {};
  \node[opt point] at (\optcut,.5) {};
  \node[opt point] at (\optbottom,0) {};
  \draw[opt guide] (\opttop,1)--(\opttop,1.14);
  \node[above=3pt] at (\opttop,1.14) {$1-\varepsilon-\tfrac12\tan\varepsilon$};
  \draw[opt guide] (\optbottom,0)--(\optbottom,-.12);
  \node[below=3pt] at (\optbottom,-.12) {$1-\varepsilon+\tfrac12\tan\varepsilon$};
  \node at (.28,.39) {$T_0^\varepsilon=T_1^\varepsilon=0$};
  \draw[opt leader] ({\optcut-.27*\opttan},.88)--(1.15,.99);
  \node[opt point] at ({\optcut-.27*\opttan},.88) {};
  \node[anchor=west,align=left] at (1.17,.99)
    {$T_0^\varepsilon=0$\\[3pt]$T_1^\varepsilon=(\cos\varepsilon,\sin\varepsilon)$};
  \draw[opt leader] (.88,.55)--(1.15,.55);
  \node[opt point] at (.88,.55) {};
  \node[anchor=west,align=left] at (1.17,.55)
    {$T_0^\varepsilon=e_1$\\[3pt]$T_1^\varepsilon=(\cos\varepsilon,\sin\varepsilon)$};
  \draw[opt leader] ({\optcut+.27*\opttan},.12)--(1.15,.15);
  \node[opt point] at ({\optcut+.27*\opttan},.12) {};
  \node[anchor=west,align=left] at (1.17,.15)
    {$T_0^\varepsilon=e_1$\\[3pt]$T_1^\varepsilon=0$};
\end{scope}
\draw[opt zero cut] (.05,-1.75)--(.7,-1.75);
\node[anchor=west] at (.95,-1.75) {$x_1=1-\varepsilon$};
\draw[opt cut] (.05,-2.4)--(.7,-2.4);
\node[anchor=west] at (.95,-2.4)
  {$\gamma^\varepsilon(x) :=x_1+\tan(\varepsilon)x_2-1+\varepsilon-\tfrac12\tan\varepsilon=0$};
\path let \p1=(current bounding box.south west),
          \p2=(current bounding box.north east),
          \n1={max(3.85cm-\x1,\x2-3.85cm)} in
  coordinate (optCenteredSW) at ({3.85cm-\n1},\y1)
  coordinate (optCenteredNE) at ({3.85cm+\n1},\y2);

\pgfresetboundingbox
\path[use as bounding box]
  (optCenteredSW) rectangle (optCenteredNE);
\end{tikzpicture}}
  \caption{Values of the optimal maps $T_0^\varepsilon, T_1^\varepsilon$ on the unit square. The dashed vertical line
  bounds $L_0^\varepsilon$, and the solid oblique line bounds $L_1^\varepsilon$.
  In the two hatched triangles, exactly one of the maps vanishes.}
  \label{fig:optimality-transport}
\end{figure}
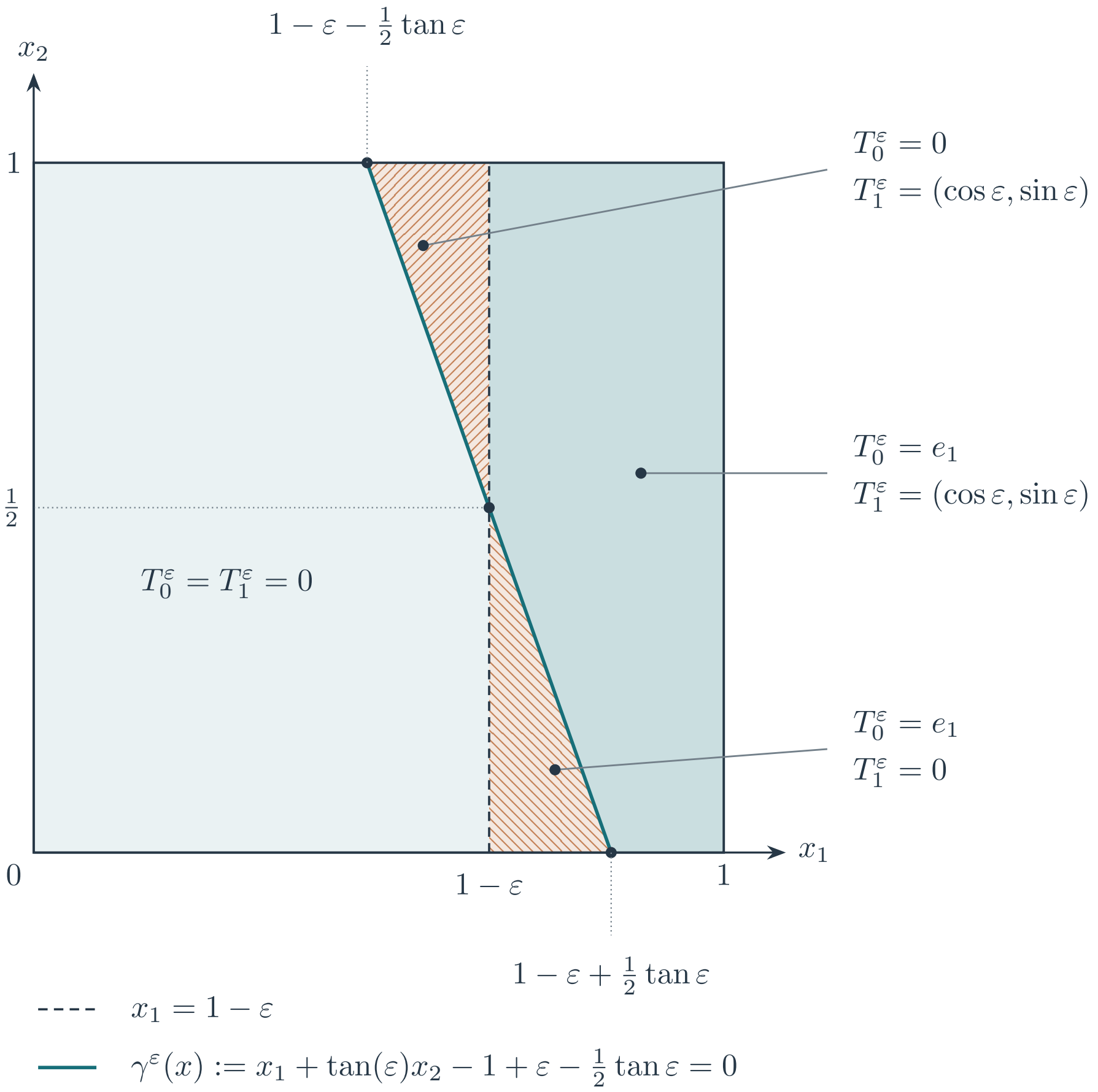

We now show that the exponent $1/4$ in Theorem~\ref{th:stab_classique} is sharp by a simple two-dimensional example. Let $\mathscr{L}^2$ denote the two-dimensional Lebesgue measure, $\rho:=\mathscr{L}^2 \resmes [0,1]^2$ be uniform on the unit square,  and  for a small parameter $\varepsilon>0$, let 
\[\rho_0^{\varepsilon}:=(1-\varepsilon) \delta_{0}+\varepsilon \delta_{e_1}, \; \rho_1^{\varepsilon}:=(1-\varepsilon) \delta_{0}+\varepsilon \delta_{(\cos(\varepsilon), \; \sin(\varepsilon))}.\]
Defining
\[\mathbb{P}_i^\varepsilon:=\frac{1}{2} \delta_\rho+ \frac{1}{2} \delta_{\rho_i^\varepsilon}, \; i=0, 1,\]
it is easy to check that 
\begin{equation}\label{W1epscarre}
\cW_1(\mathbb{P}_0^\varepsilon, \mathbb{P}_1^\varepsilon) \leq \frac{1}{2}W_1(\rho_0^{\varepsilon}, \rho_1^{\varepsilon}) \leq \frac{\varepsilon}{2} \vert e_1-(\cos(\varepsilon), \sin(\varepsilon))\vert \leq \frac{ \sqrt{2} \; \varepsilon^2}{2}.
\end{equation}
Consider the affine function $\gamma^\varepsilon$ 
\[\gamma^\varepsilon(x_1, x_2):=x_1+ \tan(\varepsilon) x_2-1+\varepsilon-\frac{\tan(\varepsilon)}{2}\]
and the convex potentials $\varphi_0^\varepsilon$ and $\varphi_1^\varepsilon$:
\[\varphi_0^\varepsilon(x_1, x_2):=\max(0, x_1 -1+\varepsilon), \; \varphi_1^\varepsilon:=\cos(\varepsilon)\max(0, \gamma^\varepsilon)\]
it is elementary to check that $T_i^\varepsilon:=\nabla \varphi_i^\varepsilon$ transports $\rho$ to $\rho_i^\varepsilon$. Hence, if we set
\[L_0^\varepsilon:=[0,1-\varepsilon]\times[0,1], \; R_0^\varepsilon:=[0,1]^2 \setminus L_0^\varepsilon\]
and 
\[L_1^\varepsilon:=[0,1]^2 \cap \{\gamma^\varepsilon \leq 0\}, \; R_1^\varepsilon:=[0,1]^2 \setminus L_1^\varepsilon,\]
we see that the optimal transport $T_0^\varepsilon$ (resp. $T_1^\varepsilon$) from $\rho$ to $\rho_0^\varepsilon$ (resp. $\rho_1^\varepsilon$) consists in sending  $L_0^\varepsilon$ (resp. $L_1^\varepsilon$) to $0$ and $R_0^\varepsilon$ (resp. $R_1^\varepsilon$) to $e_1$ (resp. $(\cos(\varepsilon), \sin(\varepsilon))$). These sets as well as the coordinates of their vertices are represented in Figure~\ref{fig:optimality-transport}.

The barycenters $\mu_0^\varepsilon$ and $\mu_1^\varepsilon$ of $\mathbb{P}_0^\varepsilon$ and $\mathbb{P}_1^\varepsilon$ coincide with the McCann's interpolants
\[\mu_i^\varepsilon:= (\frac{1}{2} \mathrm{id})_\# \nu_i^\varepsilon, \; \nu_i^\varepsilon:=(\mathrm{id}+ T_i^\varepsilon)_\# \rho.\]
One directly infers from the explicit form of the optimal transports $T_i^\varepsilon$  that  $\nu_i^\varepsilon$ is uniform on the (disjoint) union of $L_i^\varepsilon$ and some translation of $R_i^\varepsilon$, more precisely (see Figures~\ref{fig:optimality-nu-zero} and~\ref{fig:optimality-nu-one})
\[\nu_0^\varepsilon=\mathscr{L}^2 \resmes {L_0^\varepsilon\cup (e_1+R_0^\varepsilon)}, \;  \nu_1^\varepsilon=\mathscr{L}^2 \resmes {L_1^\varepsilon\cup ((\cos(\varepsilon), \sin(\varepsilon))+R_1^\varepsilon)}.\]

\begin{figure}[!t]
\centering
  \resizebox{.95\linewidth}{!}{
\begin{tikzpicture}[x=1cm,y=1cm,font=\small,text=optInk]
\path[use as bounding box](-.55,-1.4) rectangle (16.0,6.65);
\fill[optTeal!16] (0,0) rectangle (5,5.2);
\fill[optTeal!16] (8.2,0) rectangle (12.2,5.2);
\draw[opt edge] (0,0) rectangle (5,5.2);
\draw[opt edge] (8.2,0) rectangle (12.2,5.2);
\fill[optOrange!12] (8.2,0) rectangle (12.2,2);
\fill[pattern=north east lines,pattern color=optOrange!75]
  (8.2,0) rectangle (12.2,2);
\draw[optOrange,line width=.8pt] (8.2,0) rectangle (12.2,2);
\draw[opt neighborhood] (7.65,-.55) rectangle (12.75,2.55);
\draw[opt axis] (0,0)--(14.6,0) node[right] {$x_1$};
\draw[opt axis] (0,0)--(0,6.1) node[above] {$x_2$};
\OptAxisBreak{5.75}
\node[below left] at (0,0) {$0$};
\node[left] at (0,5.2) {$1$};
\node[below=4pt] at (5,0) {$1-\varepsilon$};
\node[below=4pt] at (8.2,0) {$2-\varepsilon$};
\node[below=4pt] at (12.2,0) {$2$};
\node at (2.5,2.8) {$L_0^\varepsilon$};
\node[above=5pt] at (10.2,5.2) {$e_1+R_0^\varepsilon$};
\node[fill=white,inner sep=3pt,text=optOrange] at (10.2,1.05) {$Q^\varepsilon$};
\node[anchor=west,text=optOrange] at (13.25,1.1) {$Q_\delta^\varepsilon$};
\draw[draw=optOrange,line width=.45pt] (12.75,1.1)--(13.1,1.1);
\draw[opt guide] (12.75,2.55)--(13.5,2.55);
\draw[opt guide] (12.2,2)--(13.5,2);
\node[anchor=west] at (13.6,2.55) {$\tfrac14+\delta$};
\node[anchor=west] at (13.6,2) {$\tfrac14$};
\draw[opt extension] (12.2,-.3)--(12.2,-1.0);
\draw[opt extension] (12.75,-.55)--(12.75,-1.0);
\draw[opt dimension] (12.2,-.88)--(12.75,-.88)
  node[midway,below=2pt] {$\delta$};
\end{tikzpicture}}
  \caption{The density of $\nu_0^\varepsilon$ equals $1$ on the shaded support
  and $0$ elsewhere. The hatched rectangle is
  $Q^\varepsilon=[2-\varepsilon,2]\times[0,1/4]$; the dashed outline is
  $Q_\delta^\varepsilon=Q^\varepsilon+[-\delta,\delta]^2$. Schematic representation, not to scale.}
  \label{fig:optimality-nu-zero}
\par\vspace{16pt}
\centering
  \resizebox{.95\linewidth}{!}{
\begin{tikzpicture}[x=1cm,y=1cm,font=\small,text=optInk]
\path[use as bounding box](-.55,-1.8) rectangle (16.0,6.9);
\pgfmathsetmacro{\optlefttop}{5.2+(6.8-10.1)}
\fill[optTeal!16] (0,0)--(5.2,0)--(\optlefttop,5.2)--(0,5.2)--cycle;
\draw[opt edge] (0,0)--(5.2,0)--(\optlefttop,5.2)--(0,5.2)--cycle;
\fill[optTeal!16] (10.1,.65)--(11.95,.65)--(11.95,5.85)--(6.8,5.85)--cycle;
\draw[opt edge] (10.1,.65)--(11.95,.65)--(11.95,5.85)--(6.8,5.85)--cycle;
\pgfmathsetmacro{\optintersection}{10.1+(6.8-10.1)*(2.55-.65)/(5.85-.65)}
\fill[optOrange!12] (10.1,.65)--(11.95,.65)--(11.95,2.55)--(\optintersection,2.55)--cycle;
\fill[pattern=north east lines,pattern color=optOrange!80]
  (10.1,.65)--(11.95,.65)--(11.95,2.55)--(\optintersection,2.55)--cycle;
\draw[optGray,line width=.55pt,dash pattern=on 2pt off 2pt] (8.2,0) rectangle (12.2,2);
\draw[opt neighborhood] (7.65,-.55) rectangle (12.75,2.55);
\node[anchor=north west,text=optGray] at (8.25,1.96) {$Q^\varepsilon$};
\node[anchor=west,text=optOrange] at (13.25,1.25) {$Q_\delta^\varepsilon$};
\draw[draw=optOrange,line width=.45pt] (12.75,1.25)--(13.1,1.25);
\draw[opt axis] (0,0)--(14.6,0) node[right] {$x_1$};
\draw[opt axis] (0,0)--(0,6.25) node[above] {$x_2$};
\OptAxisBreak{6.15}
\node[below left] at (0,0) {$0$};
\node[left] at (0,5.2) {$1$};
\node at (1.7,2.8) {$L_1^\varepsilon$};
\node[above=5pt] at (9.4,5.85) {$(\cos\varepsilon,\sin\varepsilon)+R_1^\varepsilon$};
\node[below=4pt] at (5.2,0) {$1-\varepsilon+\tfrac12\tan\varepsilon$};
\draw[opt guide] (11.95,5.85)--(13.5,5.85);
\node[anchor=west] at (13.6,5.85) {$1+\sin\varepsilon$};
\draw[opt guide] (12.75,2.55)--(13.5,2.55);
\node[anchor=west] at (13.6,2.55) {$\tfrac14+\delta$};
\draw[opt guide] (11.95,.65)--(13.5,.65);
\node[anchor=west] at (13.6,.65) {$\sin\varepsilon$};
\node[opt point] at (\optintersection,2.55) {};
\node[opt point] at (10.1,.65) {};
\draw[opt guide] (\optintersection,2.55)--(\optintersection,-.72);
\draw[opt guide] (10.1,.65)--(10.1,-.72);
\draw[opt guide] (11.95,.65)--(11.95,0);
\draw[opt leader] (11.95,0)--(11.95,-.24)--(13.15,-.24);
\node[anchor=west] at (13.3,-.24) {$1+\cos\varepsilon$};
\node[anchor=north east,inner sep=1pt] at (\optintersection,-.82)
  {$1+\cos\varepsilon-\varepsilon+\tan\varepsilon\bigl(\tfrac14-\delta+\sin\varepsilon\bigr)$};
\node[anchor=north west,inner sep=1pt] at (10.1,-.82)
  {$1+\cos\varepsilon-\varepsilon+\tfrac12\tan\varepsilon$};
\end{tikzpicture}}
  \caption{The density of $\nu_1^\varepsilon$ equals $1$ on the shaded support
  and $0$ elsewhere. The orange hatching shows
  $Q_\delta^\varepsilon\cap((\cos\varepsilon,\sin\varepsilon)+R_1^\varepsilon)$.
  The orange and grey dashed outlines indicate $Q_\delta^\varepsilon$
  and $Q^\varepsilon$, respectively. Schematic representation, not to scale.}
  \label{fig:optimality-nu-one}
\end{figure}

Our aim now is to show that for some positive constant $C$ and small enough $\varepsilon$ we have 
\begin{equation}\label{w2insqrteps}
W_2^2(\mu_0^\varepsilon, \mu_1^\varepsilon)=\frac{1}{4}W_2^2(\nu_0^\varepsilon, \nu_1^\varepsilon) \geq C \varepsilon.
\end{equation}
If we manage to prove \eqref{w2insqrteps}, together with \eqref{W1epscarre}, this will show that the exponent $1/4$ in Theorem~\ref{th:stab_classique} cannot be improved. To prove \eqref{w2insqrteps}, we shall take advantage of the fact that the bottom part of $(\cos(\varepsilon), \sin(\varepsilon))+R_1^\varepsilon$ is roughly twice as narrow as that of $e_1+R_0^\varepsilon$ which somehow forces a large fraction of its mass to be transported far away. Let us make it more precise by considering
\[Q^\varepsilon:=[2-\varepsilon,2] \times[0, 1/4], \; Q^\varepsilon_{\delta}:=Q^\varepsilon+ [-\delta, \delta]^2\]
for $\delta\in (0,1/4)$. By construction, we have
\begin{equation}\label{nu0massQ}
    \nu_0^\varepsilon(Q^\varepsilon)=\frac{\varepsilon}{4}.
\end{equation}
We now observe that   $Q^\varepsilon_\delta$ is disjoint from $L_1^\varepsilon$ for $\varepsilon$ small enough and then
\[\nu_1^\varepsilon(Q^\varepsilon_\delta)=\mathscr{L}^2(Q^\varepsilon_{\delta}\cap ((\cos(\varepsilon), \sin(\varepsilon))+R_1^\varepsilon)).\]
An elementary computation reveals that $Q^\varepsilon_{\delta}\cap ((\cos(\varepsilon), \sin(\varepsilon))+R_1^\varepsilon)$ lies in the vertical strip 
\[0\leq 1+\cos(\varepsilon)-x_1 \leq \varepsilon-\frac{\tan(\varepsilon)}{2}+\tan(\varepsilon)(\delta+1/4-\sin(\varepsilon)) =\varepsilon(\delta+3/4+O(\varepsilon)).\]
and of course, it is also included in the horizontal strip $x_2\in [0, 1/4+\delta]$. This immediately yields 
\begin{equation}\label{nu1small}
    \nu_1^\varepsilon(Q^\varepsilon_\delta) \leq \varepsilon (3/16+\delta+\delta^2+ O(\varepsilon)).
\end{equation}
Let now $\gamma\in \Pi(\nu_0^\varepsilon, \nu_1^\varepsilon)$, since $\vert x-y\vert^2\geq \delta^2$ when $x\in Q^\varepsilon$ and $y\in \R^2\setminus Q^\varepsilon_\delta$, we have
\[\begin{split}
    \int_{\R^2\times \R^2} \vert x-y\vert^2 d \gamma  & \geq \delta^2 \gamma(Q^\varepsilon \times (\R^2\setminus Q^\varepsilon_\delta)) \\
    & \geq \delta^2(\nu_0^\varepsilon(Q^\varepsilon)-\nu_1^\varepsilon(Q^\varepsilon_\delta))\\
&    \geq \delta^2 \varepsilon (1/16-\delta -\delta^2 +O(\varepsilon))
\end{split}\]
where we used the fact that $\gamma\in \Pi(\nu_0^\varepsilon, \nu_1^\varepsilon)$ in the second line and  \eqref{nu0massQ}- \eqref{nu1small}  in the last one. Optimizing with respect to the plan $\gamma$, we thus have 
\[W_2^2(\nu_0^\varepsilon, \nu_1^\varepsilon) \geq \delta^2 \varepsilon (1/16-\delta -\delta^2 +O(\varepsilon))\]
This proves \eqref{w2insqrteps}.


\section{Stability of Entropic Barycenters}\label{sec:entropicbar}

In this section, following \cite{Bigot19}, \cite{CEK21}, we consider a regularized variant of the Wasserstein barycenter problem. Instead of minimizing solely the Wasserstein variance, we add a penalization by the entropy of the measure $\mu$. We define the entropy of a probability measure $\mu$ on $B_R$ by: 
$$  H(\mu) = \begin{cases} \int_{B_R} \log(\mu(x))\mu(x) \, dx & \text{if } \mu \in \cP(B_R)\cap L^1(B_R) \\ +\infty & \text{otherwise} \end{cases} $$ 
where we have identified an absolutely continuous measure with respect to the Lebesgue measure with its density.
For a fixed regularization parameter $\lambda > 0$, we define the following functional on $\cP(B_R)$:
$$ F_{\mathbb{P}, \lambda}(\mu) = \int_{\cP(B_R)} \frac{1}{2} W_2^2(\mu, \rho) \, d\mathbb{P}(\rho) + \lambda H(\mu).$$
The entropic barycenter associated with $\mathbb{P}$ and parameter $\lambda$ is the unique probability measure $\mu_{\mathbb{P}, \lambda}$ minimizing the functional $F_{\mathbb{P}, \lambda}$:
$$ \mu_{\mathbb{P}, \lambda} = \mathop{\mathrm{argmin}}_{\mu \in \cP(B_R)} F_{\mathbb{P}, \lambda}(\mu).$$
Let us first summarize some useful properties of entropic barycenters from \cite{CEK21}.

\begin{theorem}[Existence, uniqueness, characterization and regularity \cite{CEK21}]\label{summaryCEK}
Let $\lambda>0$ and $\mathbb{P} \in \cP(\cP(B_R))$. There exists a unique  entropic barycenter $\mu_{\mathbb{P}, \lambda}$, moreover, this measure is absolutely continuous with respect to the Lebesgue measure with $\mathrm{supp}(\mu_{\mathbb{P}, \lambda}) = B_R$ and  satisfies the following Euler-Lagrange equation, almost everywhere on $B_R$:
\begin{equation}\label{eq:ELentropic}
\lambda \log(\mu_{\mathbb{P}, \lambda}) + \int_{\cP(B_R)} \psi_{\mu_{\mathbb{P}, \lambda} \to \rho} \, d\mathbb{P}(\rho) = 0
\end{equation}
where $\psi_{\mu_{\mathbb{P}, \lambda} \to \rho}$ denotes, for $\mathbb{P}$-almost every $\rho$, a Kantorovich potential between $\mu_{\mathbb{P}, \lambda}$ and $\rho$. Moreover, if $\mathbb{P}$ satisfies condition \ref{hyp:H0} then $\mu_{\mathbb{P},\lambda}\in L^{\infty}(B_R)$ with the explicit bound:
\begin{equation}\label{eq:boundentropicbaryc}
 \|\mu_{\mathbb{P},\lambda}\|_\infty \le \frac{M_{\mathbb{P}}}{\alpha^d_{\mathbb{P}}}.
 \end{equation}
\end{theorem}

\subsection{Quantitative convergence to classical barycenters}

Thanks to Corollary \ref{cor:forte_conv_FP}, we readily obtain a quantitative control  of the Wasserstein distance between the entropic barycenters $\mu_{\mathbb{P}, \lambda}$ and the classical Wasserstein barycenter $\mu_{\mathbb{P}}$:

\begin{proposition}\label{prop:conv_entropique}
Let $\mathbb{P} \in \cP(\cP(B_R))$ satisfying \ref{hyp:H0}.
Let $\mu_{\mathbb{P}}$ be the Wasserstein barycenter associated with $\mathbb{P}$ and $\mu_{\mathbb{P},\lambda}$ be the entropic barycenter with parameter $\lambda > 0$, then
\[ W_2(\mu_{\mathbb{P}}, \mu_{\mathbb{P},\lambda}) \le \Big(\frac{C_{d,R, M_{\mathbb{P}}} \log( M_{\mathbb{P}} \vert B_R \vert \alpha_{\mathbb{P}}^{-d})}{\alpha_{\mathbb{P}}}\Big)^{1/4} \lambda^{1/4}\]
where $C_{d, R, M_{\mathbb{P}}}$ is as in \eqref{eq:formulaforC}.
\end{proposition}

\begin{proof}
Since $\mu_{\mathbb{P}}$ is a Wasserstein barycenter, there is a family $(\psi_{\mu_{\mathbb{P}} \to \rho})_{\rho \in \mathrm{supp}(\mathbb{P})}$ of Kantorovich potentials between $\mu_{\mathbb{P}}$ and $\rho \in \mathrm{supp}(\mathbb{P})$ satisfying the Euler-Lagrange equation \eqref{eq:eulerlagrangebaryc}. Hence, by Corollary \ref{cor:forte_conv_FP}, we have
$$ F_{\mathbb{P}}(\mu_{\mathbb{P},\lambda}) - F_{\mathbb{P}}(\mu_{\mathbb{P}}) \ge \frac{\alpha_{\mathbb{P}}} {C_{d, R, M_{\mathbb{P}}}} W_2^4(\mu_{\mathbb{P},\lambda}, \mu_{\mathbb{P}}).$$
On the other hand, since $\mu_{\mathbb{P},\lambda}$ minimizes $F_{\mathbb{P},\lambda}$, we have:
$$ F_{\mathbb{P},\lambda}(\mu_{\mathbb{P}}) - F_{\mathbb{P},\lambda}(\mu_{\mathbb{P},\lambda}) \ge 0.$$

Noting that $F_{\mathbb{P},\lambda} = F_{\mathbb{P}} + \lambda H$,  and summing the two previous inequalities, we learn that
$$ \frac{\alpha_{\mathbb{P}}} {C_{d, R, M_{\mathbb{P}}}} W_2^4(\mu_{\mathbb{P},\lambda}, \mu_{\mathbb{P}}) \le \lambda \Big( H(\mu_{\mathbb{P}}) - H(\mu_{\mathbb{P},\lambda}) \Big).$$

Now, since $\mathbb{P}$ satisfies assumption \ref{hyp:H0}, by \eqref{eq:boundclassiccbaryc}, we know that $\| \mu_{\mathbb{P}} \|_\infty \le \frac{M_{\mathbb{P}}}{\alpha^d_{\mathbb{P}}}$. Therefore $H(\mu_{\mathbb{P}}) \le \log\left(\frac{M_{\mathbb{P}}}{\alpha^d_{\mathbb{P}}}\right)$. On the other hand,  $H(\mu_{\mathbb{P}, \lambda}) \ge -\log(|B_R|)$ for every $\mu \in \cP(B_R)$. Putting everything together gives
$$ W_2(\mu_{\mathbb{P},\lambda}, \mu_{\mathbb{P}}) \le \Big(\frac{C_{d,R, M_{\mathbb{P}}} \log( M_{\mathbb{P}} \vert B_R \vert \alpha_{\mathbb{P}}^{-d})}{\alpha_{\mathbb{P}}}\Big)^{1/4} \lambda^{1/4}.$$
\end{proof}

\subsection{Quantitative $\W_p$, $H(\cdot \| \cdot)$ and $L^1$ stability of entropic barycenters}

The aim of this paragraph is to gather quantitative stability results for entropic barycenters with respect to the Wasserstein distance between $\mathbb{P}$ and  $\mathbb{Q}$. We shall distinguish controls which are independent of $\lambda$ but require assumption \ref{hyp:H0} on $\mathbb{P}$ and those, stronger (in the sense that they enforce strong convergence in a quantitative manner), which rely on the regularizing effect of the entropic term and therefore depend on the regularization parameter $\lambda$.

Let us start by explaining why, under assumption \ref{hyp:H0},   the stability estimate of Theorem \ref{th:stab_classique} extends to entropic barycenters:

\begin{theorem}\label{th:stab_entropique_indep}
Let $\mathbb{P}, \mathbb{Q} \in \cP(\cP(B_R))$, $\lambda > 0$ and let $\mu_{\mathbb{P},\lambda}, \mu_{\mathbb{Q},\lambda}$ be the entropic barycenters associated with $(\mathbb{P},\lambda)$ and $(\mathbb{Q},\lambda)$. If $\mathbb{P}$ satisfies assumption \ref{hyp:H0}, then: 
\begin{equation}\label{eq:stabW2barycentropic}
 W_2(\mu_{\mathbb{P},\lambda}, \mu_{\mathbb{Q},\lambda}) \le C_{d,R, \alpha_{\mathbb{P}}, M_{\mathbb{P}}} \, \cW_1(\mathbb{P}, \mathbb{Q})^{1/4}
 \end{equation}
where $C_{d,R, \alpha_{\mathbb{P}}, M_{\mathbb{P}}}$ is the same constant as in Theorem \ref{th:stab_classique}.
\end{theorem}

\begin{proof}
The proof is similar to that of Theorem \ref{th:stab_classique}, and uses the strong convexity of $F_{\mathbb{P}, \lambda}$.
As recalled in Theorem \ref{summaryCEK}, there exists  a family of Kantorovich potentials $(\psi_{\mu_{\mathbb{P},\lambda} \to \rho})_{\rho}$ between $\mu_{\mathbb{P}, \lambda}$ and $\rho$, for $\mathbb{P}$-almost every $\rho \in \cP(B_R)$, such that the following  Euler-Lagrange holds:
$$ \int_{\cP(B_R)} \psi_{\mu_{\mathbb{P},\lambda} \to \rho}(\cdot) \, d\mathbb{P}(\rho) + \lambda \log(\mu_{\mathbb{P},\lambda}) = 0.$$
Thus for any $\mu \in \cP(B_R)$ with finite entropy, integrating the preceding equality against $(\mu - \mu_{\mathbb{P},\lambda})$ gives:
$$ \begin{aligned} 
&F_{\mathbb{P},\lambda}(\mu) - F_{\mathbb{P},\lambda}(\mu_{\mathbb{P},\lambda}) \\
&= \int_{\cP(B_R)} \left[ \frac{1}{2} W_2^2(\mu, \rho) - \frac{1}{2} W_2^2(\mu_{\mathbb{P},\lambda}, \rho) - \langle \psi_{\mu_{\mathbb{P},\lambda} \to \rho}, \mu - \mu_{\mathbb{P},\lambda} \rangle \right] d\mathbb{P}(\rho) \\
&\quad + \lambda \int_{B_R} \Big( \log(\mu(x))\mu(x) - \log(\mu_{\mathbb{P},\lambda}(x))\mu_{\mathbb{P},\lambda}(x) - [\log(\mu_{\mathbb{P},\lambda}(x)) + 1](\mu - \mu_{\mathbb{P},\lambda})(x) \Big) \, dx \\
&\ge \frac{\alpha_{\mathbb{P}}}{ C_{d, R, M_{\mathbb{P}}}} W_2^4(\mu_{\mathbb{P},\lambda}, \mu),
\end{aligned} $$
 the last inequality follows from Corollary \ref{cor:forte_conv_FP}  applied to the nonnegative integrand of the first integral over $S \subset \cP(B_R)$ and from the non-negativity (by convexity of $H$) of the second integral. In particular, taking $\mu=\mu_{\mathbb{Q},\lambda}$, we get
$$ F_{\mathbb{P},\lambda}(\mu_{\mathbb{Q},\lambda}) - F_{\mathbb{P},\lambda}(\mu_{\mathbb{P},\lambda}) \ge \frac{\alpha_{\mathbb{P}}}{ C_{d, R, M_{\mathbb{P}}}} W_2^4(\mu_{\mathbb{P},\lambda}, \mu_{\mathbb{Q},\lambda}).$$
On the other hand, the optimality  of $\mu_{\mathbb{Q},\lambda}$ for $F_{\mathbb{Q},\lambda}$ yields
$$ F_{\mathbb{Q},\lambda}(\mu_{\mathbb{P},\lambda}) - F_{\mathbb{Q},\lambda}(\mu_{\mathbb{Q},\lambda}) \ge 0.$$
Summing the two previous inequalities, the entropy terms cancel out so that, following mutatis mutandis the end of the proof of Theorem  \ref{th:stab_classique}, we arrive at the desired conclusion (with exactly the same constant $C_{d,R, \alpha_{\mathbb{P}}, M_{\mathbb{P}}})$ as in Theorem \ref{th:stab_classique}).
\end{proof}

If we now allow the constant in front of $\cW_1(\mathbb{P}, \mathbb{Q})$ to depend on the parameter $\lambda$, we are able to obtain other stronger controls taking advantage of the strict convexity of the entropy $H$. Recall the definition of the relative entropy (or Kullback-Leibler divergence) of a probability measure $\mu$ with respect to a reference probability measure $\mu_0$: 
\begin{equation}\label{Entropie_rel}
    H(\mu \| \mu_0 ) = \begin{cases} \int_{B_R} \log(\frac{d\mu}{d\mu_0}(x)) \, d\mu(x) & \text{if } \mu \ll \mu_0 \\ +\infty & \text{otherwise}. \end{cases}
\end{equation}
The following result generalizes the stability of entropic barycenters in relative entropy, already observed in \cite{Bigot19} (in the special case where $\mathbb{P}$ and $\mathbb{Q}$ are empirical measures):

\begin{proposition}\label{th:stab_entropique_H}
Let $\mathbb{P}, \mathbb{Q} \in \cP(\cP(B_R))$ and $\lambda > 0$. Consider $\mu_{\mathbb{P},\lambda}, \mu_{\mathbb{Q},\lambda}$ the entropic barycenters associated with $(\mathbb{P},\lambda)$ and $(\mathbb{Q},\lambda)$. Then:
\begin{equation}\label{eq:stabKLentropicbaryc} H(\mu_{\mathbb{Q},\lambda} \| \mu_{\mathbb{P},\lambda}) \le \frac{4R}{\lambda} \cW_1(\mathbb{P}, \mathbb{Q}).
\end{equation}
\end{proposition}

\begin{proof}
Recall that relative entropy is the Bregman divergence associated with $H$ i.e. for $\mu, \mu_0$ probabilities on $B_R$ with finite entropy, we have:
\begin{equation}\label{eq:KLisBregman}
     H(\mu \| \mu_0) = H(\mu) - H(\mu_0) - \int_{B_R} \log(\mu_0)(\mu - \mu_0).
 \end{equation}
Let $(\psi_{\mu_{\mathbb{P},\lambda} \to \rho})_{\rho}$ be a family of Kantorovich potentials between $\mu_{\mathbb{P},\lambda}$ and $\rho$ (for $\mathbb{P}$-almost every $\rho \in \cP(B_R)$) such that the Euler-Lagrange equation \eqref{eq:ELentropic} is satisfied. Then, for $\mu \in \cP(B_R)$ with finite entropy, we have:
$$ \begin{aligned} 
F_{\mathbb{P},\lambda}(\mu) - F_{\mathbb{P},\lambda}(\mu_{\mathbb{P},\lambda}) &= \int_{\cP(B_R)} \left[ \frac{1}{2} W_2^2(\mu, \rho) - \frac{1}{2} W_2^2(\mu_{\mathbb{P},\lambda}, \rho) - \langle \psi_{\mu_{\mathbb{P},\lambda} \to \rho}, \mu - \mu_{\mathbb{P},\lambda} \rangle \right] d\mathbb{P}(\rho) \\
&\quad + \lambda \left( H(\mu) - H(\mu_{\mathbb{P},\lambda}) - \int_{B_R} \log(\mu_{\mathbb{P},\lambda})(\mu - \mu_{\mathbb{P},\lambda}) \right) \\
&\ge \lambda H(\mu \| \mu_{\mathbb{P},\lambda})
\end{aligned}, $$
where the last inequality follows from the convexity of $W_2^2(\cdot, \rho)$ and \eqref{eq:KLisBregman}.
We therefore obtain:
$$ F_{\mathbb{P},\lambda}(\mu_{\mathbb{Q},\lambda}) - F_{\mathbb{P},\lambda}(\mu_{\mathbb{P},\lambda}) \ge \lambda H(\mu_{\mathbb{Q},\lambda} \| \mu_{\mathbb{P},\lambda}).$$
Adding  to the left hand side of the previous inequality, the nonnegative quantity $F_{\mathbb{Q},\lambda}(\mu_{\mathbb{P},\lambda}) - F_{\mathbb{Q},\lambda}(\mu_{\mathbb{Q},\lambda})$ and arguing exactly as in the end of the proof of Theorem \ref{th:stab_classique}, we arrive at
$$ \lambda H(\mu_{\mathbb{Q},\lambda} \| \mu_{\mathbb{P},\lambda}) \le 4R  \cW_1(\mathbb{P}, \mathbb{Q}).$$

\end{proof}
Note that under the assumptions of Proposition \ref{th:stab_entropique_H}, combining  \eqref{eq:stabKLentropicbaryc} with Pinsker's inequality, we get an $L^1$ control:
\begin{equation}\label{eq:stabL1entropicbaryc}
\|\mu_{\mathbb{Q},\lambda} - \mu_{\mathbb{P},\lambda}\|_{L^1(B_R)} \leq \frac{ 2 \sqrt{2} R^{1/2} \cW_1(\mathbb{P}, \mathbb{Q})^{1/2}}{\lambda^{1/2}}.
\end{equation}
In turn, by \eqref{eq:stabL1entropicbaryc} and the obvious inequality $W_1(\mu_{\mathbb{Q},\lambda}, \mu_{\mathbb{P},\lambda}) \leq 2R \|\mu_{\mathbb{Q},\lambda} - \mu_{\mathbb{P},\lambda}\|_{L^1(B_R)}$, we obtain
\begin{equation}\label{eq:cheapw1control}
    W_1(\mu_{\mathbb{Q},\lambda}, \mu_{\mathbb{P},\lambda}) \leq \frac{ 4 \sqrt 2 R^{3/2}  \cW_1(\mathbb{P}, \mathbb{Q})^{1/2}}{\lambda^{1/2}}.
\end{equation}
Note that this inequality has a better exponent than \eqref{eq:stabW2barycentropic} in $\cW_1(\mathbb{P}, \mathbb{Q})$ but this is of course at the expense of a negative power of $\lambda$.

\subsection{Higher-order stability of entropic barycenters} Our aim now is to obtain stability in stronger distances (enforcing strong convergence of gradients) for entropic barycenters. This somehow gives a quantitative justification for their use as robust surrogates of classical barycenters. Of course, the estimates of this paragraph will depend on $\lambda$ and will explode (in a negative power-like way) when  $\lambda$ vanishes. 

 The first-order condition \eqref{eq:ELentropic} relates the log density of $\mu_{\mathbb{P},\lambda}$ to the average of the Kantorovich potentials  $\psi_{\mu_{\mathbb{P},\lambda} \to \rho}$.  Each of these potentials is a $c$-transform for the quadratic cost; in particular, it is $2R$-Lipschitz, hence differentiable almost everywhere on $B_R$  and therefore  $\log(\mu_{\mathbb{P},\lambda})$ inherits these properties. 
Noting that since $\mu_{\mathbb{P},\lambda}$ is absolutely continuous, $\nabla \psi_{\mu_{\mathbb{P},\lambda} \to \rho}=\mathrm{id}- T_{\mu_{\mathbb{P},\lambda} \to \rho}$ where $T_{\mu_{\mathbb{P},\lambda} \to \rho}$ is Brenier's optimal transport between $\mu_{\mathbb{P},\lambda}$ and $\rho$, the differentiation of \eqref{eq:ELentropic} is legitimate (see \cite{CEK21} for details) and yields
\begin{equation}\label{eq:diffELentropic}
\lambda \nabla \log(\mu_{\mathbb{P}, \lambda}) + \int_{\cP(B_R)} (\mathrm{id}- T_{\mu_{\mathbb{P}, \lambda} \to \rho}) \, d\mathbb{P}(\rho) = 0.
\end{equation}

Let us now return to the quantitative stability of entropic barycenters. Condition  \eqref{eq:diffELentropic} allows us to bridge the stability of the gradient of the log-density to the stability  of optimal transport maps. This enables us first to state a quantitative stability results in terms of relative Fisher information:

\begin{theorem}\label{th:stab_log}
Let $\mathbb{P}, \mathbb{Q} \in \cP(\cP(B_R))$ with $\mathbb{P}$ satisfying assumption \ref{hyp:H0}. Let $\mu_{\mathbb{P},\lambda}, \mu_{\mathbb{Q},\lambda}$ be the associated entropic barycenters. Then there exists $C = C_{d, R, \alpha_{\mathbb{P}}, M_{\mathbb{P}}} > 0$ such that:
$$ \| \nabla \log(\mu_{\mathbb{P},\lambda}) - \nabla \log(\mu_{\mathbb{Q},\lambda}) \|_{L^2(\mu_{\mathbb{P},\lambda})} \le \frac{C}{\lambda} \left( \cW_1(\mathbb{P}, \mathbb{Q})^{1/4} + \cW_1(\mathbb{P}, \mathbb{Q})^{1/16} \right).$$
\end{theorem}

The quantity $\nabla \log(\mu)$ seems at first glance the most natural to consider in view of \eqref{eq:diffELentropic}; nevertheless, the relation 
$$ \nabla \log \mu = \frac{2\nabla \sqrt{\mu}}{\sqrt{\mu}} $$
suggests to consider a weighted $H^1$ norm of $\sqrt{\mu}$, because we already have a strong control of the quantity $\| \sqrt{\mu_{\mathbb{P},\lambda}} - \sqrt{\mu_{\mathbb{Q},\lambda}} \|_{L^2}$ by $\| \mu_{\mathbb{P},\lambda} - \mu_{\mathbb{Q},\lambda}\|_{L^1}$, itself controlled by inequality \eqref{eq:stabL1entropicbaryc}.

For $f \in C^{0,1}(B_R)$ and $\mu \in \cP(B_R)\cap L^1(B_R)$, we define $\| \cdot \|_{H^1(\mu)}$  by:
$$ \| f \|_{H^1(\mu)} = \| f \|_{L^2(\mu)} + \| \nabla f \|_{L^2(\mu)}.$$

We then have the following $H^1$ stability for the square roots of the entropic barycenters:

\begin{theorem}\label{th:stab_racine}
Let $\mathbb{P}, \mathbb{Q} \in \cP(\cP(B_R))$, with $\mathbb{P}$ satisfying assumption \ref{hyp:H0}. Let $0 < \lambda \le 1$ and $\mu_{\mathbb{P},\lambda}, \mu_{\mathbb{Q},\lambda}$ be the associated entropic barycenters. Then there exists $C = C_{d, R, \alpha_{\mathbb{P}}, M_{\mathbb{P}}} > 0$ such that:
$$ \| \sqrt{\mu_{\mathbb{P},\lambda}} - \sqrt{\mu_{\mathbb{Q},\lambda}} \|_{H^1(\mu_{\mathbb{P},\lambda})} \le C \left( \frac{\cW_1(\mathbb{P}, \mathbb{Q})^{1/4}}{\lambda^{5/4}} + \frac{\cW_1(\mathbb{P},\mathbb{Q})^{1/8}}{\lambda^{9/8}} \right).$$
\end{theorem}

The proof of  Theorems \ref{th:stab_log} and \ref{th:stab_racine} heavily relies on the following estimate.

\begin{lemma}\label{lem:stab_plan_transport}
Let $\mathbb{P}, \mathbb{Q} \in \cP(\cP(B_R))$ with $\mathbb{P}$ satisfying assumption \ref{hyp:H0}. Denote by $\mu_{\mathbb{P},\lambda}$ and $\mu_{\mathbb{Q},\lambda}$ the associated entropic barycenters. Then, there exists a constant $C = C_{d, R, \alpha_{\mathbb{P}}, M_{\mathbb{P}}} > 0 $ such that
$$ \left\| \int_{\cP(B_R)} T_{\mu_{\mathbb{P},\lambda} \to \rho}(\cdot) \, d\mathbb{P}(\rho) - \int_{\cP(B_R)} T_{\mu_{\mathbb{Q},\lambda} \to \rho}(\cdot) \, d\mathbb{Q}(\rho) \right\|_{L^2(\mu_{\mathbb{P},\lambda})} \le C \left( \cW_1(\mathbb{P}, \mathbb{Q})^{1/4} + \cW_1(\mathbb{P}, \mathbb{Q})^{1/16} \right).$$
\end{lemma}

\begin{proof}[Proof of Lemma \ref{lem:stab_plan_transport}]
Let us denote the entropic barycenters by $\mu_\lambda = \mu_{\mathbb{P},\lambda}$ and $\nu_\lambda = \mu_{\mathbb{Q},\lambda}$. We have:
\begin{align*}
\| \int_{\cP(B_R)} T_{\mu_\lambda \to \rho} \, d\mathbb{P}(\rho) - \int_{\cP(B_R)} T_{\nu_\lambda \to \tilde{\rho}} \, d\mathbb{Q}(\tilde{\rho}) \|^4_{L^2(\mu_\lambda)} &= \| \int_{\cP(B_R)^2} (T_{\mu_\lambda \to \rho} - T_{\nu_\lambda \to \tilde{\rho}}) \, d\gamma(\rho, \tilde{\rho}) \|^4_{L^2(\mu_\lambda)} \\
\leq \int_{\cP(B_R)^2} \| T_{\mu_\lambda \to \rho} - T_{\nu_\lambda \to \tilde{\rho}} \|^4_{L^2(\mu_\lambda)} \, d\gamma(\rho, \tilde{\rho}) 
\end{align*}
where $\gamma \in \Pi(\mathbb{P}, \mathbb{Q})$ is a transport plan between $\mathbb{P}$ and $\mathbb{Q}$.
By condition \ref{hyp:H0} on $\mathbb{P}$, we know from \eqref{eq:boundentropicbaryc} that $\|\mu_\lambda\|_\infty \le \frac{M_{\mathbb{P}}}{\alpha^d_{\mathbb{P}}}$; we can thus apply the stability Theorem for optimal transports (when both source and target measures vary) stated in \cite[Theorem 3.3]{Mer26}: 
$$\| T_{\mu_\lambda \to \rho} - T_{\nu_\lambda \to \tilde{\rho}} \|^4_{L^2(\mu_\lambda)} \leq C_1(d, R, \alpha_{\mathbb{P}}, M_{\mathbb{P}}) [ W_1(\mu_\lambda, \nu_\lambda) + W_1(\rho, \tilde{\rho})]$$
where $C_1(d, R, \alpha_{\mathbb{P}},  M_{\mathbb{P}} )=4R C_{d, R, M_{\mathbb{P}}/\alpha_{\mathbb{P}}^d}$ (with $C_{d,R,M}$ being defined by \eqref{eq:formulaforC}).

Integrating with respect to the transport plan $\gamma$, we obtain:
\begin{align*}
\| \int_{\cP(B_R)} T_{\mu_\lambda \to \rho} \, d\mathbb{P}(\rho) - \int_{\cP(B_R)} T_{\nu_\lambda \to \tilde{\rho}} \, d\mathbb{Q}(\tilde{\rho}) \|^4_{L^2(\mu_\lambda)} \\
\leq C_1(d, R, \alpha_{\mathbb{P}}, M_{\mathbb{P}} ) [ W_1(\mu_\lambda, \nu_\lambda) + \int_{\cP(B_R)^2} W_1(\rho, \tilde{\rho})\, d\gamma(\rho, \tilde{\rho})].
\end{align*}
To bound $W_1(\mu_\lambda, \nu_\lambda)$, we use  the fact that $W_1\leq W_2$ and inequality \eqref{eq:stabW2barycentropic} on the $W_2$ stability of entropic barycenters independently of the parameter $\lambda$:
\begin{equation}\label{eq:boundW1entropicbaryc1}
    W_1(\mu_\lambda, \nu_\lambda) \le C_{d,R, \alpha_{\mathbb{P}}, M_{\mathbb{P}}} \cW_1(\mathbb{P}, \mathbb{Q})^{1/4}.
\end{equation}
Minimizing over transport plans $\gamma \in \Pi(\mathbb{P}, \mathbb{Q})$, we get:
$$ \| \int_{\cP(B_R)} T_{\mu_\lambda \to \rho} \, d\mathbb{P}(\rho) - \int_{\cP(B_R)} T_{\nu_\lambda \to \tilde{\rho}} \, d\mathbb{Q}(\tilde{\rho}) \|^4_{L^2(\mu_\lambda)} \leq  C_2(d, R,  \alpha_{\mathbb{P}, M_{\mathbb{P}}}) ( \cW_1(\mathbb{P}, \mathbb{Q})^\frac{1}{4} + \cW_1(\mathbb{P}, \mathbb{Q}))$$
with the constant $ C_2(d, R, \alpha_{\mathbb{P}}, M_{\mathbb{P}})= C_1(d, R, \alpha_{\mathbb{P}}, M_{\mathbb{P}} )\max(1, C_{d,R, \alpha_{\mathbb{P}}, M_{\mathbb{P}}})$.
Finally, using the inequality $(a+b)^\frac{1}{4} \leq a^\frac{1}{4} + b^\frac{1}{4}$,
we obtain the desired inequality 
$$ \left\| \int_{\cP(B_R)} T_{\mu_{\mathbb{P},\lambda} \to \rho}(\cdot) \, d\mathbb{P}(\rho) - \int_{\cP(B_R)} T_{\mu_{\mathbb{Q},\lambda} \to \rho}(\cdot) \, d\mathbb{Q}(\rho) \right\|_{L^2(\mu_{\mathbb{P},\lambda})} \le C \left( \cW_1(\mathbb{P}, \mathbb{Q})^{1/4} + \cW_1(\mathbb{P}, \mathbb{Q})^{1/16} \right)$$
with the constant $C=C_2(d, R, \alpha_{\mathbb{P}}, M_{\mathbb{P}})^{1/4}.$
\end{proof}

Let us now prove Theorems \ref{th:stab_log} and \ref{th:stab_racine}.

\begin{proof}[Proof of Theorem \ref{th:stab_log}]
Recall that from \eqref{eq:diffELentropic}, denoting by $T_{\mu_{\mathbb{P},\lambda} \to \rho}$ the optimal transport from  $\mu_{\mathbb{P}, \lambda}$ to $\rho$, we have
$$ \lambda \nabla \log(\mu_{\mathbb{P},\lambda}) + \int_{\cP(B_R)} (\mathrm{id} - T_{\mu_{\mathbb{P},\lambda} \to \rho}) \, d\mathbb{P}(\rho) = 0.$$
Doing the same with $\mu_{\mathbb{Q},\lambda}$ and taking the difference of the two equalities, we get:
$$ \lambda \left( \nabla \log(\mu_{\mathbb{P},\lambda}) - \nabla \log(\mu_{\mathbb{Q},\lambda}) \right) =\int_{\cP(B_R)} T_{\mu_{\mathbb{P},\lambda} \to \rho} \, d\mathbb{P}(\rho)  -\int_{\cP(B_R)} T_{\mu_{\mathbb{Q},\lambda} \to \tilde{\rho}} \, d\mathbb{Q}(\tilde{\rho}).$$
Taking the $L^2(\mu_{\mathbb{P},\lambda})$ norm of both sides,  and applying Lemma \ref{lem:stab_plan_transport} yields 
$$ \| \nabla \log \mu_{\mathbb{P},\lambda} - \nabla \log \mu_{\mathbb{Q},\lambda} \|_{L^2(\mu_{\mathbb{P},\lambda})} \le \frac{C}{\lambda} \left( \cW_1(\mathbb{P}, \mathbb{Q})^{1/4} + \cW_1(\mathbb{P}, \mathbb{Q})^{1/16} \right) $$
where $C = C_{d, R, \alpha_{\mathbb{P}}, M_{\mathbb{P}}}$ is the same constant as in Lemma \ref{lem:stab_plan_transport}.

\end{proof}

\begin{rem}
    In the  proof of Lemma \ref{lem:stab_plan_transport}, instead of using \eqref{eq:stabW2barycentropic} to bound $W_1(\mu_{\mathbb{P},\lambda}, \mu_{\mathbb{Q},\lambda})$, we could have used \eqref{eq:cheapw1control} so as to obtain
    \begin{equation}\label{eq:alternativeTTcontrol}
    \left\| \int_{\cP(B_R)} T_{\mu_{\mathbb{P},\lambda} \to \rho} \, d\mathbb{P}(\rho) - \int_{\cP(B_R)} T_{\mu_{\mathbb{Q},\lambda} \to \rho} \, d\mathbb{Q}(\rho) \right\|_{L^2(\mu_{\mathbb{P},\lambda})} \le C \left( \cW_1(\mathbb{P}, \mathbb{Q})^{1/4} + \frac{\cW_1(\mathbb{P}, \mathbb{Q})^{1/8}}{\lambda^{1/8}} \right)
    \end{equation}
    and 
    \[ \| \nabla \log(\mu_{\mathbb{P},\lambda}) - \nabla \log(\mu_{\mathbb{Q},\lambda}) \|_{L^2(\mu_{\mathbb{P},\lambda})} \le \frac{C}{\lambda} \left( \cW_1(\mathbb{P}, \mathbb{Q})^{1/4} + \frac{\cW_1(\mathbb{P}, \mathbb{Q})^{1/8}}{\lambda^{1/8}} \right).\]
    \end{rem}

\begin{proof}[Proof of Theorem \ref{th:stab_racine}]
Let us simply set $\mu = \mu_{\mathbb{P},\lambda}$ and $\nu = \mu_{\mathbb{Q},\lambda}$. By \eqref{eq:diffELentropic}, we have
$$ \lambda \nabla \log \mu + \int_{\cP(B_R)} (\mathrm{id} - T_{\mu \to \rho}) \, d\mathbb{P}(\rho) = 0.$$
Multiplying by $\sqrt{\mu}$ and using the identity $\sqrt{\mu} \nabla \log \mu = 2 \nabla \sqrt{\mu}$, we obtain
$$ 2\lambda \nabla \sqrt{\mu} + \sqrt{\mu} \int_{\cP(B_R)} (\mathrm{id} - T_{\mu \to \rho}) \, d\mathbb{P}(\rho) = 0.$$
Likewise
$$ 2\lambda \nabla \sqrt{\nu} + \sqrt{\nu} \int_{\cP(B_R)} (\mathrm{id} - T_{\nu \to \tilde{\rho}}) \, d\mathbb{Q}(\tilde{\rho}) = 0,$$
and taking the difference yields:
\begin{align*}
2\lambda( \nabla \sqrt{\mu} - \nabla \sqrt{\nu} ) &= \sqrt{\nu} \int_{\cP(B_R)} (\mathrm{id} - T_{\nu \to \tilde{\rho}}) \, d\mathbb{Q}(\tilde{\rho}) - \sqrt{\mu} \int_{\cP(B_R)} (\mathrm{id} - T_{\mu \to \rho}) \, d\mathbb{P}(\rho) \\
 &= \underbrace{(\sqrt{\nu} - \sqrt{\mu}) \int_{\cP(B_R)} (\mathrm{id} - T_{\nu \to \tilde{\rho}}) \, d\mathbb{Q}(\tilde{\rho})}_{:=I}\\
& +\underbrace{ \sqrt{\mu} \big( \int_{\cP(B_R)}T_{\mu \to \rho} \, d\mathbb{P}(\rho) - \int_{\cP(B_R)} T_{\nu \to \tilde{\rho}} \, d\mathbb{Q}(\tilde{\rho}) \big)}_{:=II}.
\end{align*}
Hence
\begin{equation}\label{eq:boundwitIandII}
\| \sqrt{\mu} - \sqrt{\nu} \|_{H^1(\mu)}  \leq  \| \sqrt{\mu} - \sqrt{\nu} \|_{L^2(\mu)} + \frac{1}{2\lambda}\Vert I \Vert_{L^2(\mu)}+ \frac{1}{2\lambda} \Vert II \Vert_{L^2(\mu)}.
\end{equation}
Thanks to \eqref{eq:boundentropicbaryc} and the inequality  $|\sqrt{b} - \sqrt{a}|^2 \le |b - a|$, we firstly have
\begin{align*}
\| \sqrt{\mu} - \sqrt{\nu} \|_{L^2(\mu)}^2 &\le \| \mu - \nu \|_{L^1(\mu)} \\
&\le \frac{M_{\mathbb{P}}}{\alpha_{\mathbb{P}}^d} \| \mu - \nu \|_{L^1} \\
&\le \frac{ 2  \sqrt{2} M_{\mathbb{P}} R^{1/2} \cW_1(\mathbb{P}, \mathbb{Q})^{1/2}}{\alpha_{\mathbb{P}}^d\lambda^{1/2}}
\end{align*}
where we used \eqref{eq:stabL1entropicbaryc} in the last line.
Secondly, since $\|\mathrm{id}-T_{\nu \to \tilde{\rho}}\|_\infty \le 2R$ for all $\tilde{\rho} \in \cP(B_R)$, we also have
$$ \|I\|_{L^2(\mu)} \le 2 R  \| \sqrt{\mu} - \sqrt{\nu} \|_{L^2(\mu)} \leq C \frac{1}{\lambda^{1/4}} \cW_1(\mathbb{P}, \mathbb{Q})^{1/4} $$
where, from now, on $C=C_{d,  R, \alpha_{\mathbb{P}}, M_{\mathbb{P}}}$ denotes a positive constant (possibly changing from line to line) depending on $d, R, \alpha_{\mathbb{P}},$ and $M_{\mathbb{P}}$ but not on $\lambda$. 
Finally, using again \eqref{eq:boundentropicbaryc} and \eqref{eq:alternativeTTcontrol}, we learn that
$$ \|II\|_{L^2(\mu)} \le C \left( \cW_1(\mathbb{P}, \mathbb{Q})^{1/4} + \frac{\cW_1(\mathbb{P}, \mathbb{Q})^{1/8}}{\lambda^{1/8}} \right).$$
With \eqref{eq:boundwitIandII}, we therefore get 
\begin{align*}
\| \sqrt{\mu} - \sqrt{\nu} \|_{H^1(\mu)} 
&\le C \left( \frac{1}{\lambda^{1/4}} \cW_1(\mathbb{P}, \mathbb{Q})^{1/4} + \frac{1}{\lambda^{5/4}} \cW_1(\mathbb{P}, \mathbb{Q})^{1/4} + \frac{1}{\lambda} \left( \cW_1(\mathbb{P}, \mathbb{Q})^{1/4} + \frac{\cW_1(\mathbb{P}, \mathbb{Q})^{1/8}}{\lambda^{1/8}} \right) \right) \\
&\le C \left( \frac{\cW_1(\mathbb{P}, \mathbb{Q})^{1/4}}{\lambda^{5/4}} + \frac{\cW_1(\mathbb{P},\mathbb{Q})^{1/8}}{\lambda^{9/8}} \right)
\end{align*}
provided that $0 < \lambda \le 1$.
\end{proof}


\section{Stability of doubly entropic barycenters}\label{sec:doublyentropicbar}

In this final section, we address the stability of doubly entropic barycenters introduced by Chizat \cite{Chi25}. For computational tractability reasons, \cite{Cut13}, \cite{CuturiPeyre}, it has become common in the literature  to replace the squared Wasserstein distance $W_2^2$ with its entropically regularized version $\operatorname{OT}_\lambda$, defined for $\lambda > 0$ by:
\begin{equation}\label{eq:defOTlambda}
    \forall \mu, \rho \in \cP(B_R), \quad \operatorname{OT}_\lambda(\mu, \rho) = \inf_{\pi \in \Pi(\mu, \rho)} \left\{ \int_{B_R^2} \frac{1}{2}|x-y|^2 \, d\pi(x,y) + \lambda H(\pi \| \mu \otimes \rho)\right\}
\end{equation}
where $H(\cdot \| \cdot)$ is the relative entropy functional (see \eqref{Entropie_rel}). We refer to the survey \cite{Leo14} and the lecture notes \cite{Nutz2021} for an overview of entropic optimal transport.  Given $\mathbb{P} \in \cP(\cP(B_R))$ and $\lambda>0$, we consider the doubly entropic variance functional  
\[
    F_{\mathbb{P}, \lambda, \lambda}(\mu) = \int_{\cP(B_R)} \operatorname{OT}_\lambda(\mu, \rho) \, d\mathbb{P}(\rho) + \lambda H(\mu)
\]
and denote by $\mu_{\mathbb{P}, \lambda, \lambda}$ the doubly entropic barycenter associated with $\mathbb{P}$ and $\lambda$ i.e. the unique minimizer of $F_{\mathbb{P}, \lambda, \lambda}$ over $\cP(B_R)$. We refer to  \cite{Chi25} for motivations, properties and characterizations of these barycenters. 

\smallskip

We will first give a convergence rate for the Wasserstein distance between doubly entropic barycenters and classical barycenters in paragraph \ref{subsec:cvbmupllmup}. Our approach to the quantitative stability of doubly entropic barycenters will be based on an extension  to the entropic case of the strong convexity estimate from \cite[Theorem 2.2]{Mer26}, see  Theorems \ref{th:forte_conv_K_rho_lambda}  and \ref{th:forte_conv_OT_lambda_plus_entropie} below. To prove these estimates, we will need to recall basic facts from entropic optimal transport and study some properties of soft  Legendre transforms which will be the object of paragraph \ref{subsec:softlegendre}. Paragraph \ref{sec:forte_convexe_K_rho_lambda} will be devoted to the strong convexity of $\mathcal{K}_\rho^\lambda$, the entropic version of the Kantorovich functional. Finally quantitative stability for doubly entropic baycenters will be addressed in paragraph \ref{subsec:quantitativedoubleyentropic}.

\subsection{Quantitative convergence to classical barycenters}\label{subsec:cvbmupllmup}
Under condition \ref{hyp:H0} on $\mathbb{P}$, we have the following quantitative convergence rate of doubly entropic barycenters towards the classical Wasserstein barycenter:
\begin{proposition}\label{prop:conv_double_entropique}
Let $\mathbb{P} \in \cP(\cP(B_R))$ satisfying \ref{hyp:H0}. Denoting by $\mu_{\mathbb{P}}$ (respectively $\mu_{\mathbb{P}, \lambda, \lambda}$) the corresponding Wasserstein (resp. doubly entropic) barycenter, then there exists 
$C=C_{d, R, \alpha_{\mathbb{P}}, M_{\mathbb{P}}}$ such that
\[
    W_2(\mu_{\mathbb{P}}, \mu_{\mathbb{P}, \lambda, \lambda}) \le C \lambda^{\frac{1}{4}} (|\log \lambda| + 1)^{\frac{1}{4}}.
\]
\end{proposition}

\begin{proof}
Since $\mu_{\mathbb{P}}$ is a Wasserstein barycenter, there is a family $(\psi_{\mu_{\mathbb{P}} \to \rho})_{\rho \in \mathrm{supp}(\mathbb{P})}$ of Kantorovich potentials between $\mu_{\mathbb{P}}$ and $\rho \in \mathrm{supp}(\mathbb{P})$ satisfying the Euler-Lagrange equation \eqref{eq:eulerlagrangebaryc}. Hence, by Corollary \ref{cor:forte_conv_FP}, we have
$$ F_{\mathbb{P}}(\mu_{\mathbb{P},\lambda,\lambda}) - F_{\mathbb{P}}(\mu_{\mathbb{P}}) \ge \frac{\alpha_{\mathbb{P}}} {C_{d, R, M_{\mathbb{P}}}} W_2^4(\mu_{\mathbb{P},\lambda,\lambda}, \mu_{\mathbb{P}}).$$
On the other hand, the fact that $\mu_{\mathbb{P},\lambda, \lambda}$ minimizes $F_{\mathbb{P},\lambda,\lambda}$ gives
$$ F_{\mathbb{P},\lambda, \lambda}(\mu_{\mathbb{P}}) - F_{\mathbb{P},\lambda, \lambda}(\mu_{\mathbb{P},\lambda, \lambda}) \ge 0.$$
Summing these two inequalities, we learn that
\begin{align*}
    \frac{\alpha_{\mathbb{P}}} {C(d, R, M_{\mathbb{P}})} W_2^4(\mu_{\mathbb{P},\lambda, \lambda}, \mu_{\mathbb{P}}) &\le \int_{\cP(B_R)} \left( \frac{1}{2} W_2^2(\mu_{\mathbb{P}}, \rho) - \operatorname{OT}_\lambda(\mu_{\mathbb{P}}, \rho) \right) d\mathbb{P}(\rho) \\
    &\quad + \int_{\cP(B_R)} \left( \operatorname{OT}_\lambda(\mu_{\mathbb{P}, \lambda, \lambda}, \rho) - \frac{1}{2} W_2^2(\mu_{\mathbb{P}, \lambda, \lambda}, \rho) \right) d\mathbb{P}(\rho) \\
    &\quad + \lambda(H(\mu_{\mathbb{P}}) - H(\mu_{\mathbb{P}, \lambda, \lambda})).
\end{align*}

The first term is non-positive. By \cite[Cor. 3.14]{EN23}, see also \cite{CPT23}, we have $\operatorname{OT}_\lambda - \frac{1}{2} W_2^2 \le C_{d, R} \lambda (|\log \lambda| + 1)$. Thus the second term is bounded by the same quantity. Finally, since $\mathbb{P}$ satisfies assumption \ref{hyp:H0}, by \eqref{eq:boundclassiccbaryc}, we know that  $\| \mu_{\mathbb{P}} \|_\infty \le \frac{M_{\mathbb{P}}}{\alpha^d_{\mathbb{P}}}$. Therefore $H(\mu_{\mathbb{P}}) \le \log\left(\frac{M_{\mathbb{P}}}{\alpha^d_{\mathbb{P}}}\right)$. On the other hand,  $H(\mu) \ge -\log(|B_R|)$ for every $\mu \in \cP(B_R)$. Putting everything together gives
$$W_2(\mu_{\mathbb{P},\lambda}, \mu_{\mathbb{P}}) \le  C_{d, R, \alpha_{\mathbb{P}}, M_{\mathbb{P}}}  \lambda^{\frac{1}{4}} (|\log \lambda| + 1)^{\frac{1}{4}}.$$

\end{proof}

\subsection{Entropic optimal transport, soft Legendre transforms and Kantorovich functional}\label{subsec:softlegendre} 
It is well-known that the unique optimal entropic plan $\pi^\lambda$ between $\mu$ and $\rho$ i.e the minimizer in $\operatorname{OT}_\lambda(\mu, \rho)$ can be written in exponential form as 
\[ d\pi^{\lambda}(x,y)=\exp\Big(\frac{-c(x,y)+\psi(x)+\phi(y)}{\lambda}\Big) d \mu(x) \; d \rho(y)\]
where $c(x,y)=\frac{1}{2} \vert x-y\vert^2$ and  the potentials $\phi$ and $\psi$ (usually refered to as Schr\"odinger potentials in the literature) are such that $\pi\in \Pi(\mu, \rho)$ i.e. solve the Schr\"odinger system:
\begin{align}
    \psi(x)=-\lambda \log\Big(\int_{B_R} e^{\frac{-c(x,y)+\phi(y)}{\lambda}} d\rho(y) \Big)=:\phi^{c, \lambda, \rho}(x), \mbox{ for $\mu$-a.e. $x$,}
     \nonumber\\
     \phi(y)=-\lambda \log\Big(\int_{B_R} e^{\frac{-c(x,y)+\psi(x)}{\lambda}} d\mu(x) \Big)=:\psi^{c, \lambda, \mu}(y), \mbox{ for $\rho$-a.e. $y$.}
      \label{eq:schrosys}
    \end{align}
This system which is at the heart of the popular Sinkhorn  algorithm has been intensively studied. It turns out to be well-posed, indeed, as shown in \cite{CCL24}, extended to all $x$ and $y$ in $B_R$, \eqref{eq:schrosys} is uniquely solvable (with a Lipschitz with respect to $W_1$ dependence on the marginals $\mu$ and $\rho$) in $C(B_R)\times C(B_R)$, up to equivalence relation: $(\psi_1, \phi_1) \sim (\psi_2, \phi_2)$ if and only if there exists a constant $C$ such that $\psi_1=\psi_2+C$ and $\phi_1=\phi_2-C$. Thus, we may denote by $\psi^{\lambda}_{\mu\to \rho}$ (Schr\"odinger potential between $\rho$ and $\mu$)  and $\phi^{\lambda}_{\rho\to \mu}$ (Schr\"odinger potential between $\mu$ and $\rho$) the unique (up to the equivalence relation mentioned above) solution in $C(B_R)\times C(B_R)$ of the system 
\[\psi=\phi^{c, \lambda, \rho} \mbox{ on $B_R$}, \;  \phi=\psi^{c, \lambda, \mu} \mbox{ on $B_R$}.\]
Note also that Schr\"odinger potentials are smooth and  $2R$-Lipschitz on $B_R$.

It is also well-known that \eqref{eq:schrosys} is nothing but the Euler-Lagrange system for  the dual formulation of $\operatorname{OT}_\lambda(\mu, \rho)$
\begin{equation*}
\operatorname{OT}_\lambda(\mu, \rho)=\lambda+ \sup_{(\phi, \psi)\in C(B_R)^2} \left\{ \langle \psi, \mu \rangle+ \langle \phi, \rho \rangle-\lambda\int_{B_R^2} e^{\frac{-c(x,y)+\psi(x)+\phi(y)}{\lambda}} d\mu(x)  d\rho(y)\right\}.
    \end{equation*}
    It may also be convenient to replace this dual formulation by the semi-dual version simply obtained by optimizing with respect to $\phi$ (i.e. by replacing $\phi$ by $\psi^{c, \lambda, \mu}$):
 \begin{equation}\label{eq:semidualeotclassic}
\operatorname{OT}_\lambda(\mu, \rho)=\sup_{\psi\in C(B_R)} \left\{ \langle \psi, \mu \rangle+ \langle \psi^{c, \lambda, \mu} \, \rho \rangle\right\}. 
\end{equation}
Note that since the supremum is achieved by the $2R$-Lipschitz potential $\psi^{\lambda}_{\mu \to \rho}$, by Kantorovich-Rubinstein inequality, we readily obtain that for every $\mu\in \cP(B_R)$, $\rho\in \cP(B_R)\mapsto \operatorname{OT}_\lambda(\mu, \rho)$ is $2R$-Lipschitz for $W_1$. 
Moreover, \eqref{eq:semidualeotclassic}, together with the well-posedness of the Schr\"odinger system implies that $\operatorname{OT}_\lambda$ is differentiable (see \cite{CCL24}, \cite{Chi25} for details) in the sense that for any $\mu, \rho$ and $\alpha, \beta$ in $\cP(B_R)$, we have
\begin{align*}
\lim_{\varepsilon \downarrow 0} \varepsilon^{-1}(\operatorname{OT}_\lambda(\mu+\varepsilon(\alpha-\mu), \rho)-\operatorname{OT}_\lambda(\mu, \rho))=\int_{B_R} \psi^{\lambda}_{\mu\to \rho} d(\alpha-\mu), \\
 \lim_{\varepsilon \downarrow 0} \varepsilon^{-1}(\operatorname{OT}_\lambda(\mu, \rho+ \varepsilon(\beta-\rho))-\operatorname{OT}_\lambda(\mu, \rho))=\int_{B_R} \phi^{\lambda}_{\rho\to \mu} d(\beta-\rho)
 \end{align*}
 in the sequel, we shall simply write these conditions as
 \begin{equation}\label{eq:gradOTlambdapot}
 \nabla_\mu \operatorname{OT}_\lambda(\mu, \rho)=\psi^{\lambda}_{\mu\to \rho},\;  \nabla_\rho \operatorname{OT}_\lambda(\mu, \rho)=\phi^{\lambda}_{\rho\to \mu}.
 \end{equation}
Simply performing the change of variable $v:=\frac{1}{2}\vert \cdot\vert^2-\psi$ in \eqref{eq:semidualeotclassic} yields 
\begin{equation}\label{eq:semidualeotclassicalabrenier}
\operatorname{OT}_\lambda(\mu, \rho)=\frac{1}{2}(m_2(\mu)+m_2(\rho))-\inf_{v\in C(B_R)}\{ \langle v, \mu\rangle +\lambda \int_{B_R} \log\Big(\int_{B_R} e^{\frac{x \cdot y-v(x)}{\lambda}} d\mu(x) \Big) d \rho(y)\}.
    \end{equation}
 If, in addition, $\mu$ is absolutely continuous with respect to Lebesgue and $\log(\mu)\in C(B_R)$ (a property which is automatically satisfied by doubly entropic barycenters, as we will recall later), we can make a last change of variable $u:=v-\lambda \log(\mu)\in C(B_R)$ and obtain the following identity which is the entropic analogue of \eqref{eq:Brenierdf}:
\begin{equation}\label{Eq:Semi-dual-entropique}
    \operatorname{OT}_\lambda(\mu, \rho) + \lambda H(\mu) = \frac{1}{2}(m_2(\mu) + m_2(\rho)) - \min_{u \in C(B_R)} \left\{ \langle u, \mu \rangle + \mathcal{K}_{\rho}^{\lambda}(u) \right\}
\end{equation}
where we define for every $u\in C(B_R)$, the soft Legendre transform of $u$, $u^{*, \lambda}$ by
\[
    u^{*, \lambda}(y) := \lambda \log \left( \int_{B_R} e^{\frac{1}{\lambda}(x \cdot y - u(x))} \, dx \right), \forall y\in \R^d.
\]
and the entropic Kantorovich functional $\mathcal{K}_{\rho}^{\lambda}$  by: 
\[
    \mathcal{K}_{\rho}^{\lambda}(u) := \int_{B_R} u^{*, \lambda}(y) \, d\rho(y), \;  \forall u\in C(B_R).
\]
Note (again assuming that $\log(\mu)\in C(B_R)$) that the infimum in \eqref{Eq:Semi-dual-entropique} is achieved by $u=\frac{1}{2}\vert \cdot\vert^2-\psi^{\lambda}_{\mu\to \rho}-\lambda \log(\mu)$. Recalling \eqref{eq:gradOTlambdapot} this rewrites as 
\begin{equation}\label{eq:gradient_OT+H}
    \nabla_\mu (\operatorname{OT}_{\lambda}(\cdot, \rho)+\lambda H)(\mu)=\frac{1}{2}\vert \cdot\vert^2-u.  
\end{equation}

Given $u\in C(B_R)$,  and $y\in \R^d$ let us now define the  probability (density) $p_u^y$ on $B_R$ by
\begin{equation}\label{eq:defdepuy}
    p_u^y(x) := \frac{e^{\frac{1}{\lambda}(x \cdot y - u(x))}}{\int_{B_R} e^{\frac{1}{\lambda}(z \cdot y - u(z))} \, dz}, \; x\in B_R.
\end{equation}
By construction $u^{*, \lambda}$ is smooth and its gradient and hessian are respectively given by the barycenter and $\lambda^{-1}$ times the variance of $p_u^y$ (so that $u^{*, \lambda}$ is strictly convex). Note that $p_u^y$ can be rewritten as $p_u^y(x) = e^{\frac{1}{\lambda}\left(x \cdot y - u(x) - u^{*, \lambda}(y)\right)}$. If, in addition, $u$ is an optimal potential in (\ref{Eq:Semi-dual-entropique}) then $u^{*, \lambda}=\frac{1}{2}\vert \cdot \vert^2-\phi^{\lambda}_{\rho \to \mu}$  and $p_u^y$ is the density of the disintegration $\pi^\lambda(\cdot | y)$ of the entropic transport plan $\pi^\lambda$ with respect to $\rho$. We then define $T^\lambda=T^{\lambda}_{\rho \to \mu}$ the entropic analogue of the Brenier's  map between $\rho$ and $\mu$ by
\begin{equation}\label{eq:entropic_brenier_map}
    T^\lambda(y) := \int_{B_R} z \, d\pi^\lambda(z|y), \; y \in \R^d,
\end{equation}
and observe that $T^\lambda(y) = \int_{B_R} z \, dp_u^y(z) = \nabla u^{*,\lambda}(y)$ so that, as Brenier's map, $T^\lambda$ has a convex potential. Moreover, it follows from \cite{KLM25}, that the soft Kantorovich functional is indeed smooth and in particular satisfies
\begin{equation}\label{eq:formulagradkrholambda}
    \nabla \mathcal{K}_{\rho}^{\lambda}(u) =-\int_{B_R} p_u^y \; d\rho(y) 
    \end{equation}
in the sense that for every $v\in C(B_R)$
\[\lim_{t\to 0} t^{-1} (\mathcal{K}_{\rho}^{\lambda}(u+tv)-\mathcal{K}_{\rho}^{\lambda}(u))=-\int_{B_R} \int_{B_R} v(x) p_u^y(x)dx d\rho(y).\] 
We will also need the following result concerning the (usual) Legendre transform $(u^{*, \lambda})^*$ of the soft Legendre transform $u^{*, \lambda}$ i.e.
\[(u^{*, \lambda})^*(x):=\sup_{y\in \R^d} \{x\cdot y- u^{*, \lambda}(y)\}, \; x\in \R^d.\]

\begin{lemma}\label{lem:legendre_soft_max}
For $u \in C(B_R)$, the Legendre transform of $u^{*, \lambda}$ admits the following representation
\begin{equation}\label{eq:softenvelopeconvexe}
    (u^{*, \lambda})^*(x) = \inf_{\nu \in \cP(B_R)} \left\{ \langle u, \nu \rangle + \lambda H(\nu) : \int_{B_R} z \, d\nu(z) = x \right\}
\end{equation}
in particular, $(u^{*, \lambda})^*(x)=+\infty$  if $\vert x\vert \ge R$. If $\vert x \vert <R$, the infimum in \eqref{eq:softenvelopeconvexe} is attained at $\nu^* = p_u^{y^*}$ where $y^*\in \R^d$ is the unique root of the equation $\nabla u^{*, \lambda}(y^*) = x$.
\end{lemma}

\begin{proof}
Let us recall that the Donsker-Varadhan variational representation of relative entropy says that for $f \in C(B_R)$ and $\mu_0 \in \cP(B_R)$, one has: 
\[
\log \left( \int_{B_R} e^{f(x)} d\mu_0(x) \right) = \sup_{\nu \in \cP(B_R)} \left\{ \int_{B_R} f(z) \, d\nu(z) - H(\nu \| \mu_0) \right\}.
\]
Taking  $f(z) = \frac{y\cdot z - u(z)}{\lambda}$ and $\mu_0$ the normalized Lebesgue measure on $B_R$, we thus have: 
\begin{align*}
    u^{*, \lambda}(y) &= \sup_{\nu \in \cP(B_R)} \left\{ \int_{B_R} \big( \langle y, z \rangle - u(z) \big) \, d\nu(z) - \lambda H(\nu ) \right\}\\
    &=\sup_{x\in B_R} \left\{ y\cdot x -\inf_{\nu \in \cP(B_R)} \left\{ \langle u, \nu \rangle + \lambda H(\nu) : \int_{B_R} z \, d\nu(z) = x \right\}\right\}
\end{align*}
which means that  $u^{*, \lambda}$ is the Legendre transform of
\[x\mapsto g(x):=\inf_{\nu \in \cP(B_R)} \left\{ \langle u, \nu \rangle + \lambda H(\nu) : \int_{B_R} z \, d\nu(z) = x \right\}\]
since $g$ is convex and lsc, by Fenchel-Moreau Theorem, we have $(u^{*, \lambda})^*=g^{**}=g$ which proves formula \eqref{eq:softenvelopeconvexe} as well as the fact that $(u^{*, \lambda})^*(x)=+\infty$ when $\vert x\vert\geq R$ (if $\vert x\vert>R$ this is clear since no probability concentrated on $B_R$ has $x$ as center of mass and when $\vert x\vert=R$ only the Dirac mass at $x$, which is of infinite entropy has this property). 

Now let $\delta\in (0,R)$ and $x$ be such that $\vert x \vert\leq R-\delta$. Observe that for $y\in \R^d\setminus\{0\}$ $K_\delta(y):=\{z\in B_R\; : \; z\cdot y \geq (R-\delta/2) \vert y\vert\}$ has a positive Lebesgue measure $m_\delta$ (independent of $y$), so that
\[u^{*, \lambda}(y) \geq \lambda \log(m_\delta)+(R-\delta/2)\vert y\vert-\Vert u\Vert_{\infty}\]
hence
\[ x\cdot y-u^{*, \lambda}(y) \leq -\frac{\delta}{2} \vert y\vert-\lambda \log(m_\delta)+ \Vert u\Vert_{\infty} \to -\infty \mbox{ as $\vert y \vert\to \infty$}\]
so that $y\mapsto x\cdot y-u^{*, \lambda}(y)$, has a (unique by strict convexity of $u^{*, \lambda}$) maximizer which we denote by $y^*$ and is characterized by the optimality condition $x=\nabla u^{*, \lambda}(y^*)=\int_{B_R} z d \nu^*(z)$ for $\nu^*:=p_u^{y^*}$. It remains to check that $\nu^*$ (which has barycenter $x$) is optimal in \eqref{eq:softenvelopeconvexe}, to see this, we proceed by direct evaluation of the functional at $\nu^*$:
\begin{align*}
    \int_{B_R} &u(z) \, d\nu^*(z) + \lambda H(\nu^*)  
    = \int_{B_R} u(z) \, d\nu^*(z) + \lambda \int_{B_R} \log\left( \frac{\exp\left(\frac{y^* \cdot z  - u(z)}{\lambda}\right)}{\int_{B_R} \exp\left(\frac{y^* \cdot w  - u(w)}{\lambda}\right) dw} \right) d\nu^*(z) \\
    &= \int_{B_R} u(z) \, d\nu^*(z) + \int_{B_R} \Big( y^* \cdot z  - u(z) \Big) d\nu^*(z) - \lambda \log \int_{B_R} \exp\left(\frac{y^* \cdot w  - u(w)}{\lambda}\right) dw. \\
    &= y^* \cdot x  - u^{*, \lambda}(y^*)=(u^{*, \lambda})^*(x)
\end{align*}
and conclude by formula \eqref{eq:softenvelopeconvexe}.

\end{proof}

\subsection{Strong convexity of $\mathcal{K}_\rho^\lambda$}
\label{sec:forte_convexe_K_rho_lambda}

The following result generalizes the quantitative strong convexity  of $\mathcal{K}_{\rho}$ from  \cite[Theorem 2.2]{Mer26} to the entropic case.

\begin{theorem}\label{th:forte_conv_K_rho_lambda}
Let $\rho \in \cP(B_R) \cap L^\infty(B_R)$. For all $u, u_0 \in C(B_R)$,
\[
   B_{\mathcal{K}_\rho^\lambda}(u, u_0):= \mathcal{K}_{\rho}^{\lambda}(u) - \mathcal{K}_{\rho}^{\lambda}(u_0) - \langle \nabla \mathcal{K}_{\rho}^{\lambda}(u_0), u - u_0 \rangle \ge \frac{\|\nabla u^{*, \lambda} - \nabla u_0^{*, \lambda}\|_{L^2(\rho)}^4 } {C_{d, R, \|\rho\|_\infty} }
\]
 where $C_{d, R,M}$ is as in \eqref{eq:formulaforC}.
\end{theorem}

\begin{proof}
By \eqref{eq:formulagradkrholambda}, we have:
\[
    B_{\mathcal{K}_\rho^\lambda}(u, u_0) = \int_{B_R} \left( u^{*, \lambda}(y) - u_0^{*, \lambda}(y) + \langle u - u_0, p_{u_0}^y  \rangle \right) d\rho(y).
\]
Let $T_0^\lambda = \nabla u_0^{*, \lambda}$, we then have 
\[
    u_0^{*, \lambda}(y) + (u_0^{*, \lambda})^*(T_0^\lambda(y)) = y \cdot T_0^\lambda(y), \; \forall y\in \R^d.
\]
On the other hand, since by construction $T_0^\lambda(y)$ is the barycenter of $p_{u_0}^y$, it follows from Lemma \ref{lem:legendre_soft_max} that
\[    (u_0^{*, \lambda})^*(T_0^\lambda(y)) = \langle u_0, p_{u_0}^y \rangle + \lambda H(p_{u_0}^y )\]
and thanks to \eqref{eq:softenvelopeconvexe}, we infer 
\[
    (u^{*, \lambda})^*(T_0^\lambda(y)) \le \langle u, p_{u_0}^y \rangle + \lambda H(p_{u_0}^y).
\]
Thus, we can  bound from below the integrand in $B_{\mathcal{K}_\rho^\lambda}(u, u_0)$ as follows
\[\begin{split}
    u^{*, \lambda}(y) - u_0^{*, \lambda}(y) + \langle u - u_0, p_{u_0}^y \rangle 
    &= u^{*, \lambda}(y) + (u_0^{*, \lambda})^*(T_0^\lambda(y)) - y \cdot T_0^\lambda(y) + \langle u - u_0, p_{u_0}^y \rangle \\
    &= u^{*, \lambda}(y) + \lambda H(p_{u_0}^y ) - y \cdot T_0^\lambda(y) + \langle u, p_{u_0}^y \rangle \\
    &\ge u^{*, \lambda}(y) + (u^{*, \lambda})^*(T_0^\lambda(y)) - y \cdot T_0^\lambda(y).
\end{split}\]
Therefore:
\[  B_{\mathcal{K}_\rho^\lambda}(u, u_0) \ge \int_{B_R} \left( u^{*, \lambda}(y) + (u^{*, \lambda})^*(T_0^\lambda(y)) - y \cdot T_0^\lambda(y) \right) d\rho(y).\]
At this stage, we notice that the right hand side of the previous inequality is exactly the integral of the \emph{Fenchel duality gap} considered in \cite[Lemma A.1]{Mer26}, so exactly the same proof  as in Appendix A of \cite{Mer26}, gives the desired inequality
\[
   C_{d, R, \|\rho\|_\infty}  B_{\mathcal{K}_\rho^\lambda}(u, u_0) \ge \| \nabla u^{*, \lambda} - T_0^\lambda \|^4_{L^2(\rho)}.
\]
\end{proof}

 The Bregman divergence of $\operatorname{OT}_\lambda(\cdot, \rho) + \lambda H$  between $\mu$ and $\mu_0$ in $\cP(B_R)$ is by definition given by:
\begin{equation*}
    B_{\operatorname{OT}_\lambda(\cdot, \rho) + \lambda H }(\mu, \mu_0) := \operatorname{OT}_\lambda(\mu, \rho) + \lambda H(\mu) - \operatorname{OT}_\lambda(\mu_0, \rho) - \lambda H(\mu_0) - \langle \psi_{\mu_0 \to \rho}^\lambda + \lambda\log(\mu_0), \mu - \mu_0 \rangle.
\end{equation*}
Thanks to the duality relation \eqref{Eq:Semi-dual-entropique}, we can reformulate Theorem \ref{th:forte_conv_K_rho_lambda} as a strong convexity statement for $\operatorname{OT}_\lambda(\cdot, \rho) + \lambda H$:
\begin{corollary}\label{cor:forte_conv_OT_lambda}
Let $\rho \in \cP(B_R) \cap L^\infty(B_R)$. For all $\mu, \mu_0 \in \cP(B_R)\cap L^1(B_R)$  such that $\log \mu$, $\log \mu_0$ $\in C(B_R)$, we have
\begin{equation}
    B_{\operatorname{OT}_\lambda(\cdot, \rho) + \lambda H }(\mu, \mu_0) \ge \frac{ \|T^\lambda - T_0^\lambda \|_{L^2(\rho)}^4}{C_{d, R, \|\rho\|_\infty}} 
\end{equation}
where $T^\lambda$, $T_0^\lambda$ are  the entropic Brenier's maps associated with $\rho$ and $\mu, \mu_0$, respectively (defined in (\ref{eq:entropic_brenier_map})) and $C_{d, R, \|\rho\|_\infty}$ is as in \eqref{eq:formulaforC}.
\end{corollary}

\begin{proof}
Let $u \in C(B_R)$ (resp. $u_0 \in C(B_R)$) be optimal potentials in the duality formula (\ref{Eq:Semi-dual-entropique}) for the pair $(\rho, \mu)$ (resp. $(\rho, \mu_0)$). Recall first that, by (\ref{eq:gradient_OT+H}), we have
\[
    \nabla_\mu (\operatorname{OT}_\lambda(\cdot, \rho) + \lambda H)(\mu_0) = \frac{1}{2}|\cdot|^2 - u_0
\]
and since $u$ in a minimizer in \eqref{Eq:Semi-dual-entropique}, it satisfies the optimality  condition $\nabla \mathcal{K}_\rho^\lambda(u) = - \mu$. Thus,  using the duality formula (\ref{Eq:Semi-dual-entropique}), we get:
\begin{align*}
    B_{\operatorname{OT}_\lambda(\cdot, \rho) + \lambda H}(\mu, \mu_0) &= \left( \frac{1}{2} m_2(\mu) + \frac{1}{2} m_2(\rho) - \langle u, \mu \rangle - \mathcal{K}_\rho^\lambda(u) \right) \\
    &\quad - \left( \frac{1}{2} m_2(\mu_0) + \frac{1}{2} m_2(\rho) - \langle u_0, \mu_0 \rangle - \mathcal{K}_\rho^\lambda(u_0) \right) \\
    &\quad - \left\langle \frac{1}{2}|\cdot|^2 - u_0, \mu - \mu_0 \right\rangle \\
    &= \mathcal{K}_\rho^\lambda(u_0) - \mathcal{K}_\rho^\lambda(u) + \langle u_0 - u, \mu \rangle \\
    &= B_{\mathcal{K}_\rho^\lambda}(u_0, u) \\
    &\ge \frac{ \|T^\lambda - T_0^\lambda\|_{L^2(\rho)}^4}{C_{d, R, \|\rho\|_\infty}},
\end{align*}
where the last inequality follows from Theorem \ref{th:forte_conv_K_rho_lambda} since $\nabla u^{*,\lambda} = T^\lambda$ and $\nabla u_0^{*,\lambda} = T_0^\lambda$.
\end{proof}

Our goal now is to control the Wasserstein distance between $\mu$ and $\mu_0$ by $B_{\operatorname{OT}_\lambda(\cdot, \rho) + \lambda H }(\mu, \mu_0)$ and an error term vanishing as $\lambda\to 0$. For this, we need a control of the distance between Brenier's maps and their entropic counterparts, as in Theorem 4.3 in \cite{Mer26}.

\begin{lemma}\label{lem:raffinement_brenier_entropic_mer26}
Let $\rho, \mu \in \cP(B_R)$, with $\rho \in L^\infty(B_R)$ and $\mu$ with finite entropy. Let $T:=T_{\rho\to \mu}$ and $T^{\lambda}:=T^{\lambda}_{\rho\to \mu}$ be  respectively Brenier's map and Brenier's entropic map  between $\rho$ and $\mu$, then, there exists $C=C_{d, R, \|\rho\|_\infty}>0$ such that
\[
    \| T^\lambda - T \|_{L^2(\rho)}^4 \le C \lambda (H(\rho) + H(\mu) +2 \log(\vert B_R\vert)+1).
\]
\end{lemma}

\begin{proof}
Note first that if $\lambda\geq 1$, the desired inequality is obvious with the constant $C=(2R)^4$. We therefore  from now on assume that $\lambda\in (0,1)$. Introducing  $P_\lambda = \int_{B_R^2} \frac{1}{2}|x-y|^2 \, d\pi^\lambda(x,y) - \frac{1}{2} W_2^2(\mu, \rho) \ge 0$, thanks to Theorem 4.3 of \cite{Mer26} and Jensen's inequality, there is a positive constant $C=C_{d, R, \|\rho\|_\infty}$ such that
 \begin{equation}\label{eq:ineq4.5-0}
    \| T^\lambda - T \|_{L^2(\rho)}^4 \le C P_\lambda.
\end{equation} 
To control $P_\lambda$, we will use some results from \cite{MS25} and \cite{EN23}. In the remainder of this proof, $C_{d,R}$ will denote a positive constant only depending on $d$ and $R$ but whose value may change from line to line. 
Firstly,  \cite[Cor. 3.14]{EN23} yields
\begin{equation*}
   P_\lambda+ \lambda  H(\pi^\lambda \| \rho \otimes \mu)= \operatorname{OT}_\lambda(\rho, \mu) - \frac{1}{2} W_2^2(\rho, \mu) \le- \frac{d}{2} \lambda \log \lambda + C_{d,R} \lambda
\end{equation*}
and since  $H(\pi^\lambda \|  \rho \otimes \mu) = H(\pi^\lambda) - H(\rho) - H(\mu)$, we get
\begin{equation}\label{eq:ineq4.5-1}
    P_\lambda + \lambda H(\pi^\lambda) \le -\frac{d}{2} \lambda \log \lambda + \lambda (H(\rho) + H(\mu) + C_{d,R}).
\end{equation}
Secondly, since the support of $\pi^\lambda$ is included in $B_R\times B_R$,  Proposition 3.5 of \cite{MS25} allows us to bound  $H(\pi^\lambda)$ from below using $P_\lambda$ as
\begin{equation}\label{eq:ineq4.5-2}
    H(\pi^\lambda) \ge - \frac{d}{2} \log(P_\lambda) -\log(\vert B_R\vert)-\frac{d}{2} \log\Big( \frac{4 \pi e}{d}\Big).
\end{equation}
Combining \eqref{eq:ineq4.5-1} and \eqref{eq:ineq4.5-2},  we obtain

\[
    \frac{P_\lambda}{\lambda} - \frac{d}{2} \log\left( \frac{P_\lambda}{\lambda} \right) \le H(\rho) + H(\mu) + C_{d,R},
\]
so, thanks to the concavity inequality $\log(t) \le \log(d) + \frac{1}{d}( t - d)$, we  arrive at 
\[
    \frac{P_\lambda}{ 2\lambda} - \frac{d}{2}\log(d) + \frac{d}{2} \le H(\rho) + H(\mu) + C_{d,R}.
\]
Therefore:
\[
    P_\lambda \le C_{d,R} \lambda \left( H(\rho) + H(\mu) + 2 \log(\vert B_R\vert)+1 \right)
\]
and we conclude by \eqref{eq:ineq4.5-0}.

\end{proof}

Thanks to Corollary \ref{cor:forte_conv_OT_lambda} and Lemma \ref{lem:raffinement_brenier_entropic_mer26}, we now have: 

\begin{theorem}\label{th:forte_conv_OT_lambda_plus_entropie}
Let $\rho \in \cP(B_R) \cap L^\infty(B_R)$.  There exists a positive constant  $C=C_{d, R, \|\rho\|_\infty}>0$ such that, for all $\mu, \mu_0 \in \cP(B_R)\cap L^1(B_R)$ satisfying $\log \mu$, $\log \mu_0$ $\in C(B_R)$, one has
\[W_2(\mu, \mu_0)^4 \leq C \Big(
B_{\operatorname{OT}_\lambda(\cdot, \rho) + \lambda H(\cdot) }(\mu, \mu_0) + \lambda(H(\mu)+  H(\mu_0)+ 2 \log(\vert B_R\vert)+1)\Big).
\]
\end{theorem}

\begin{proof}
Since $(T_{\rho \to \mu}, T_{\rho \to \mu_0})_\# \rho\in \Pi(\mu, \mu_0)$ we have (denoting by $C$ a positive constant that only depends on $d$, $R$ and $\|\rho\|_\infty$):
\begin{align*}
    & W_2(\mu, \mu_0)^4 \leq \|T_{\rho \to \mu} - T_{\rho \to \mu_0}\|_{L^2(\rho)}^4 \\
&\leq 27 (\|T_{\rho \to \mu} - T^\lambda_{\rho \to \mu}\|_{L^2(\rho)}^4+ \|T^\lambda_{\rho \to \mu} - T^\lambda_{\rho \to \mu_0}\|_{L^2(\rho)}^4 + \|T^\lambda_{\rho \to \mu_0} - T_{\rho \to \mu_0}\|_{L^2(\rho)}^4 ) \\
&\leq C \Big(
B_{\operatorname{OT}_\lambda(\cdot, \rho) + \lambda H(\cdot) }(\mu, \mu_0) + \lambda(H(\mu)+  H(\mu_0)+ 2 \log(\vert B_R\vert)+1)\Big)
\end{align*}
where we used the triangle inequality and the convexity of $t\mapsto t^4$ in the second line and Corollary \ref{cor:forte_conv_OT_lambda}, Lemma \ref{lem:raffinement_brenier_entropic_mer26} and $H(\rho) \leq \log(\|\rho\|_\infty)$ in the last line.

\end{proof}

\subsection{Quantitative stability}\label{subsec:quantitativedoubleyentropic}
Theorem \ref{th:forte_conv_OT_lambda_plus_entropie} easily translates into a quantitative stability statement for doubly entropic barycenters provided $\mathbb{P}$ satisfies assumption \ref{hyp:H0}. For the sake of completeness we also include a stability statement in relative entropy (degenerating when $\lambda \rightarrow 0$).

\begin{theorem}
Let $\mathbb{P}, \mathbb{Q} \in \cP(\cP(B_R))$. Let $\mu_{\mathbb{P}, \lambda, \lambda}$ and $\mu_{\mathbb{Q}, \lambda, \lambda}$ be the doubly entropic barycenters associated with $\mathbb{P}$ and $\mathbb{Q}$. Then:
\begin{equation}\label{eq:stab_double_bar_relative_entro}
    H(\mu_{\mathbb{Q},\lambda, \lambda} \| \mu_{\mathbb{P},\lambda, \lambda}) \le \frac{4R}{\lambda} \cW_1(\mathbb{P}, \mathbb{Q}).
\end{equation}
Furthermore, if $\mathbb{P}$ satisfies \ref{hyp:H0}, there exists $C=C_{d, R, \alpha_{\mathbb{P}}, M_{\mathbb{P}}}>0$ such that 
\begin{equation}\label{eq:stab_double_bar_W_2}
    W_2(\mu_{\mathbb{P}, \lambda, \lambda}, \mu_{\mathbb{Q}, \lambda, \lambda}) \le C \left[\cW_1(\mathbb{P}, \mathbb{Q}) + \lambda (|\log \lambda|+1) \right]^{\frac{1}{4}}.
\end{equation}
\end{theorem}

\begin{proof}
To shorten notations, set $\mu^\lambda = \mu_{\mathbb{P}, \lambda, \lambda}$ and $\nu^\lambda = \mu_{\mathbb{Q}, \lambda, \lambda}$. The first-order optimality condition for doubly entropic barycenters (see \cite[Theorem 2.5]{Chi25} for a proof when $\mathbb{P}$ is discrete, which readily extends to our continuous setting) reads as 
\begin{equation}\label{eq:premier_ordre_bar_double_entropic}
\nabla F_{\mathbb{P}, \lambda, \lambda}(\mu^\lambda) := \int_{\cP(B_R)} \psi^\lambda_{ \mu^\lambda \rightarrow \rho} \; d\mathbb{P}(\rho) + \lambda \log \mu^\lambda = 0
\end{equation}
for some family of Schr\"odinger potentials $\psi^\lambda_{ \mu^\lambda \rightarrow \rho}$. The proof of  (\ref{eq:stab_double_bar_relative_entro}) then proceeds exactly as in Proposition \ref{th:stab_entropique_H}, using the convexity of $\operatorname{OT}_\lambda(\cdot, \rho)$ and 
\begin{equation}\label{eq:OTl1lipw1}
    \max\Big(\mathrm{Lip}_{W_1} (\operatorname{OT}_\lambda( \mu^\lambda, \cdot)), \mathrm{Lip}_{W_1} (\operatorname{OT}_\lambda(\nu^\lambda, \cdot ))\Big) \le 2R 
\end{equation}
as we recalled in paragraph \ref{subsec:softlegendre}. 

\smallskip

Let us now suppose that $\mathbb{P}$ satisfies \ref{hyp:H0} and prove (\ref{eq:stab_double_bar_W_2}). Note first that the optimality condition (\ref{eq:premier_ordre_bar_double_entropic}) guarantees that the log densities of $\mu^\lambda$ and $\nu^\lambda$ are continuous, then we can apply Theorem \ref{th:forte_conv_OT_lambda_plus_entropie} and get
\begin{align*}
    F_{\mathbb{P}, \lambda, \lambda}(\nu^\lambda) - F_{\mathbb{P}, \lambda, \lambda}(\mu^\lambda) 
    &=  F_{\mathbb{P}, \lambda, \lambda}(\nu^\lambda) - F_{\mathbb{P}, \lambda, \lambda}(\mu^\lambda) + \langle \nabla F_{\mathbb{P}, \lambda, \lambda}(\mu^\lambda), \nu^\lambda - \mu^\lambda \rangle \\
    &\ge \int_S B_{\operatorname{OT}_\lambda(\cdot, \rho) + \lambda H}(\nu^\lambda, \mu^\lambda) \, d\mathbb{P}(\rho) \\
    &\ge \frac{\alpha_{\mathbb{P}}}{C_{d, R, M_{\mathbb{P}}}} W_2^4(\nu^\lambda, \mu^\lambda) - \lambda (H(\nu^\lambda) + H(\mu^\lambda) + 2 \log(\vert B_R\vert) + 1) 
\end{align*}
where $S \subset L^\infty(B_R)$ is given by assumption \ref{hyp:H0} on $\mathbb{P}$.
Since $\nu^\lambda$ and $\mu^\lambda$ are doubly entropic barycenters, thanks to Lemma \ref{lem:entropie_barycentre_double}, proved below,  we have
\[
    \lambda(H(\nu^\lambda) + H(\mu^\lambda)) \le C_{d, R} \lambda (|\log \lambda| + 1).
\]
Therefore, for some positive constant $C = C_{d,R,\alpha_{\mathbb{P}}, M_{\mathbb{P}} }$, it holds that
\[
    \frac{ W_2^4(\nu^\lambda, \mu^\lambda)}{C} \leq  \lambda(|\log \lambda| + 1)+ F_{\mathbb{P}, \lambda, \lambda}(\nu^\lambda) - F_{\mathbb{P}, \lambda, \lambda}(\mu^\lambda).
\]
Adding the nonnegative quantity  $F_{\mathbb{Q}, \lambda, \lambda}(\mu^\lambda) - F_{\mathbb{Q}, \lambda, \lambda}(\nu^\lambda)$ to the right-hand side, yields

\begin{align*}
    \frac{W_2^4(\nu^\lambda, \mu^\lambda)}{C} & \leq  \lambda(|\log \lambda| + 1)  + \int_{\cP(B_R)} (\operatorname{OT}_\lambda( \nu^\lambda, \rho) - \operatorname{OT}_\lambda(\mu^\lambda, \rho)) \, d(\mathbb{P} - \mathbb{Q})(\rho) \\
    &\le \lambda(|\log \lambda| + 1) + \mathrm{Lip}_{W_1} (\operatorname{OT}_\lambda(\nu^\lambda, \cdot) - \operatorname{OT}_\lambda(\mu^\lambda, \cdot)) \cW_1(\mathbb{P}, \mathbb{Q})
    \\ 
    & \le  \lambda(|\log \lambda| + 1) + 4R \cW_1(\mathbb{P}, \mathbb{Q})
\end{align*}
where we used the Kantorovich-Rubinstein inequality in the second line  and \eqref{eq:OTl1lipw1} in the third one.

\end{proof}

In the previous proof, we have used the following a priori bound on the entropy of  doubly entropic barycenters:
\begin{lemma}\label{lem:entropie_barycentre_double}
If $\mu_{\mathbb{P}, \lambda, \lambda}$ is the doubly entropic barycenter associated with $\mathbb{P} \in \cP(\cP(B_R))$, then there exists $C=C_{d, R}>0$ such that
\[
    H(\mu_{\mathbb{P}, \lambda, \lambda}) \le C_{d, R} (|\log \lambda| + 1).
\]
\end{lemma}

\begin{proof}
Let us set $\mu^\lambda = \mu_{\mathbb{P}, \lambda, \lambda}$ and  let $\mu^0$ be a (non-regularized) Wasserstein barycenter associated with $\mathbb{P}$. For $\varepsilon \in (0,R)$, let $\mu^0_\varepsilon = \chi_\varepsilon * ({{(\mathrm{Proj}_{B_{R - \varepsilon}})}_\#\mu^0)} \in \cP(B_R)$ where $\mathrm{Proj}_{B_{R-\varepsilon}}$ is the orthogonal projection on $B_{R-\varepsilon}$ and $\chi_\varepsilon = \varepsilon^{-d}\chi(\frac{\cdot}{\varepsilon})$ is a mollifier with $\chi$ a smooth probability density on $\R^d$ with support in $B_1$. We have $\frac{1}{2}W_2^2 \le \operatorname{OT}_\lambda$, thus by minimality of $\mu^\lambda$:
\begin{align*}
    \int_{\cP(B_R)} \frac{1}{2} W_2^2(\mu^\lambda, \rho) \, d\mathbb{P}(\rho) + \lambda H(\mu^\lambda) &\le \int_{\cP(B_R)} \operatorname{OT}_\lambda(\mu^\lambda, \rho) \, d\mathbb{P}(\rho) + \lambda H(\mu^\lambda) \\
    &\le \int_{\cP(B_R)} \operatorname{OT}_\lambda(\mu^0_\varepsilon, \rho) \, d\mathbb{P}(\rho) + \lambda H(\mu^0_\varepsilon),
\end{align*}
hence
\begin{align}
    \lambda H(\mu^\lambda) &\le \int_{\cP(B_R)} \left(\operatorname{OT}_\lambda(\mu^0_\varepsilon, \rho) - \frac{1}{2} W_2^2(\mu^\lambda, \rho)\right) d\mathbb{P}(\rho) + \lambda H(\mu^0_\varepsilon) \nonumber \\
    &= \int_{\cP(B_R)} \left(\operatorname{OT}_\lambda(\mu^0_\varepsilon, \rho) - \frac{1}{2} W_2^2(\mu^0_\varepsilon, \rho)\right) d\mathbb{P}(\rho) \nonumber \\
    &\quad + \int_{\cP(B_R)} \left(\frac{1}{2} W_2^2(\mu^0_\varepsilon, \rho) - \frac{1}{2} W_2^2(\mu^0, \rho)\right) d\mathbb{P}(\rho) \nonumber\\
    &\quad + \int_{\cP(B_R)} \left(\frac{1}{2} W_2^2(\mu^0, \rho) - \frac{1}{2} W_2^2(\mu^\lambda, \rho)\right) d\mathbb{P}(\rho) + \lambda H(\mu^0_\varepsilon) \nonumber\\
    & \leq C_{d,R} \lambda (|\log \lambda| + 1)+ 4 R\varepsilon+\lambda H(\mu^0_\varepsilon) \label{eq:ineqintermentropybarlambda}
\end{align}
where  we used:
\begin{itemize}
\item the fact that, by \cite[Cor. 3.14]{EN23}, we have $\operatorname{OT}_\lambda - \frac{1}{2} W_2^2 \le C_{d, R} \lambda (|\log \lambda| + 1)$, 

\item  the triangle inequality and the fact that all our measures have support in $B_R$, so that for all $\rho \in \cP(B_R)$: 
\[
    W_2^2(\mu^0_\varepsilon, \rho) - W_2^2(\mu^0, \rho) \leq 4RW_2(\mu^0_\varepsilon,\mu^0) \leq 8R \varepsilon,
\]
\item the fact that $\mu^0$ minimizes the variance $F_{\mathbb{P}}$.
\end{itemize}

 Finally, since $h(t) = t\log(t)$ is convex, by Jensen inequality we have:
\[
    H(\mu^0_\varepsilon) \le H(\chi_\varepsilon) = d|\log(\varepsilon)| + H(\chi).
\]
Taking $\varepsilon = \min(R/2,\lambda)$ and recalling \eqref{eq:ineqintermentropybarlambda}, we finally obtain (up to changing the value of the constant $C_{d, R}$):
\[
    \lambda H(\mu^\lambda) \le C_{d, R} \lambda (|\log \lambda| + 1).
\]
\end{proof}

\bibliographystyle{plain}

\bibliography{References}

\end{document}